\documentclass[onecolumn,12pt]{IEEEtran}
\renewcommand{\baselinestretch}{1.5}

\IEEEoverridecommandlockouts
\usepackage[T1]{fontenc}

\usepackage{amsmath,amssymb}

\usepackage{url}

\usepackage{algorithm,algpseudocode}
\makeatletter
    \algnewcommand{\Statenl}[1]{\Statex \hskip\ALG@thistlm\hskip\algorithmicindent #1}
    \algnewcommand{\LComment}[1]{\Statex \hskip\ALG@thistlm\parbox[t]{\dimexpr\linewidth-\ALG@thistlm}{\strut $\triangleright$  #1\strut}}
    \algnewcommand{\wComment}[1]{\Statex \hskip\ALG@thistlm\hskip\algorithmicindent \parbox[t]{\dimexpr\linewidth-\ALG@thistlm-\algorithmicindent}{\strut $\triangleright$  #1\strut}}
\makeatother

\usepackage{subcaption}
\usepackage{color,tikz}
\definecolor{green}{rgb}{0,0.4,0}
\newcommand{\R}{\mathbb{R}}
\newcommand{\tr}{\intercal}

\newcommand{\cA}{\mathcal{A}}
\newcommand{\cB}{\mathcal{B}}
\newcommand{\cC}{\mathcal{C}}
\newcommand{\cD}{\mathcal{D}}
\newcommand{\cE}{\mathcal{E}}
\newcommand{\cG}{\mathcal{G}}
\newcommand{\cI}{\mathcal{I}}
\newcommand{\cM}{\mathcal{M}}
\newcommand{\cN}{\mathcal{N}}
\newcommand{\cO}{\mathcal{O}}
\newcommand{\cS}{\mathcal{S}}
\newcommand{\cU}{\mathcal{U}}
\newcommand{\cV}{\mathcal{V}}

\newcommand{\cX}{\mathcal{X}}

\newcommand{\bA}{\bar{A}}
\newcommand{\bB}{\bar{B}}

\newcommand{\bE}{\bar{E}}
\newcommand{\bH}{\bar{H}}

\newcommand{\bM}{\bar{M}}

\newtheorem{definition}{Definition}
\newtheorem{theorem}{Theorem}

\newtheorem{lemma}{Lemma}

\newtheorem{remark}{Remark}

\title{\huge Distributed Verification of Structural Controllability for Linear Time-Invariant Systems}

\author{%
    Jo\~ao~Carvalho~$^{\dagger,\ast}$\quad %
    S\'ergio~Pequito~$^{\dagger,\ddagger,\ast}$\quad %
    A.~Pedro~Aguiar~$^{\dagger,\diamond}$\quad \\
    Soummya~Kar~$^{\ddagger}$\quad %
    Karl~H.~Johansson~$^{\mathsection}$ %
    \thanks{%
    This work was partially supported by grant SFRH/BD/33779/2009, from Funda\c{c}\~ao para a Ci\^encia e a Tecnologia (FCT) and the CMU-Portugal (ICTI) program, and by projects CONAV/FCT-PT (PTDC/EEACRO/113820/2009), FCT (PEst-OE/EEI/LA0008/2013) and MORPH (EU FP7 No.  288704).}%
    \thanks{${}^{\ast}$%
    The first two authors made equal contributions to the research and writing of this paper.}%
    \thanks{${}^{\dagger}$%
    Institute for Systems and Robotics, Instituto Superior T\'ecnico, University of Lisbon, Lisbon, Portugal.}%
    \thanks{${}^{\ddagger}$%
    Department of Electrical and Computer Engineering, Carnegie Mellon University, Pittsburgh, PA 15213.}%
    \thanks{${}^{\diamond}$%
    Department of Electrical and Computer Engineering, Faculty of Engineering, University of Porto (FEUP), Portugal.}%
    \thanks{${}^{\mathsection}$%
ACCESS Linnaeus Center, School of Electrical Engineering, KTH Royal Institute of Technology, Stockholm, Sweden.}}

\begin{document}

\maketitle

\begin{abstract}
    Motivated by the development and deployment of large-scale dynamical systems,  often composed of geographically distributed smaller subsystems, we address the problem of verifying their controllability in a distributed manner.  In this work we study controllability in the structural system theoretic sense, \emph{structural} controllability. In other words, instead of focusing on a specific numerical system realization, we provide guarantees for equivalence classes of linear time-invariant systems on the basis of their structural sparsity patterns, i.e.,\ location of zero/nonzero entries in the plant matrices.
    To this end, we first propose several necessary and/or sufficient conditions to ensure structural controllability of the overall system, on the basis of the structural patterns of the subsystems and their interconnections.  The proposed verification criteria are shown to be efficiently implementable (i.e.,  with polynomial time complexity in the number of the state variables and inputs) in two important subclasses of interconnected dynamical systems: \emph{similar} (i.e., every subsystem has the same structure), and \emph{serial} (i.e., every subsystem outputs to at most one other subsystem).
    Secondly, we provide a distributed algorithm to verify structural controllability for  interconnected dynamical systems. The proposed distributed algorithm is efficient and implementable at the subsystem level; the algorithm is iterative, based on communication among (physically) interconnected subsystems, and requires only local model and interconnection knowledge at each subsystem.
\end{abstract}


\vspace{-.2cm}
\section{Introduction}

Recent years have witnessed a growing use of large-scale dynamical systems, notably, those with a modular structure~\cite{DecentralizedInterconnected,Davison1983169,Davison1977109}, such as content delivery networks, social networks, robot swarms, and smart grids.  Such systems, often geographically distributed, are composed of smaller subsystems, which we may refer to as \emph{agents}. A typical concern is that of ensuring that the system, as a whole, performs as intended.  More than often,  while analyzing these interconnected dynamical systems, which in this paper are considered as continuous linear-time invariant (LTI) subsystems,  the exact parameters of the plant matrices are unknown.  Therefore, we focus on the zero/nonzero pattern of the system's plant, which we refer to as \emph{sparsity pattern}, and, in particular, on structural counterpart of controllability, i.e., \emph{structural controllability}\footnote{The sparsity pattern (or \emph{structural matrix}) of a real matrix $M$ is a binary matrix $\bM$, satisfying $\bM_{i,j} = 0$ if and only if $M_{i,j} = 0$.  We then call the pair of structural matrices $(\bA,\bB)$ associated to an LTI system, a \emph{structural system}, and say that it is \emph{structurally controllable} if and only if there exists a controllable pair of real matrices $(A',B')$ with the same sparsity pattern as $(\bA,\bB)$.  In fact, we notice that if $(\bA,\bB)$ is structurally controllable, then almost all pairs of real matrices $(A,B)$ with the same sparsity pattern as $(\bA,\bB)$ are controllable~\cite{dionSurvey}.}~\cite{dionSurvey}.

It is worthwhile noting that these agents may be \emph{homogeneous} or \emph{heterogeneous}, from the structural point of view.  When the agents are homogeneous, their plants and connections have the same sparsity pattern and the system is referred to as a \emph{similar system}.  Otherwise, the agents are heterogeneous and two possible scenarios are conceivable: $(i)$ an agent may receive input from (possibly several) other agents but it only feeds to one other agent, in which case the overall system is referred to as \emph{serial}; and $(ii)$ the interactions between agents can be arbitrary, typically expected in setups where the inputs to the agents is broadcasted from the others.  In this paper,  all the above subclasses of interconnected dynamical systems are of interest and explored in detail.  More precisely,  several necessary and/or sufficient conditions are proposed to ensure key properties of the system; further, these properties can be verified by resorting to efficient algorithms, i.e., with polynomial time complexity in the number of state and input variables.

In some applications, the problem of composability is particularly relevant; for example, a swarm of robots possessing similar structure where the interaction (possibly through communication or information exchange) topology may change over time, or where robots may join or leave the swarm over time.  Then, the existence of necessary and/or sufficient conditions on the structure and interconnection between these agents contribute to \emph{controllability-by-design} schemes, i.e., we ensure that by inserting an agent into the interconnected dynamical system the dynamical system continues to be controllable.  In addition, we can also specify with which agents an agent should interact with such that those conditions hold.

A swarm of robots can also be composed of a variety of heterogeneous agents in which case controllability-by-design is also important. Yet, due to constraints on the communication range, the interaction between agents is merely local, even if some additional information is known.  Therefore, in the context of serial systems we can equip each subsystem with the capability of inferring if the entire system is structurally controllable or not. In other words,  distributed algorithms that rely only on interactions between a subsystem and its neighbors, where information about their structure may be shared.  In particular, if we equip the robots in the swarm with actuation capabilities that can be activated when the interconnected dynamical system is not structurally controllable, we can render this interconnected dynamical system structurally controllable.

Nonetheless, imposing \emph{a priori} knowledge of the structure of the interconnections in the system (for instance, whether it is a serial system) can be restrictive, so distributed algorithms that aim to verify structural controllability of general interconnected dynamical systems are desirable.  Hereafter, we provide such an algorithm: it requires the interaction between a subsystem and its neighbors, but it does not require (directly) sharing the structure of the subsystems involved.  Instead, it requires only partial information about its structure, which is also desirable from a privacy viewpoint.  The proposed scheme is also particularly suitable to other applications such as the smart grid of the future, that consists of entities described by subsystems deployed over large distances; in particular, notice that in these cases, the different entities may not be willing to share information about their structure due to security or privacy reasons.

\vspace{0.2cm}
\noindent \textbf{Related Work: }
Structural controllability was introduced by Lin~\cite{Lin_1974} in the context of single-input single-output (SISO) systems, and extended to multi-input multi-output (MIMO) systems by Shields and Pearson~\cite{Shields_Pearson:1976}.  A recent survey of the results in structural systems theory, where several necessary and sufficient conditions are presented, can be found in~\cite{dionSurvey}.

In this paper,  we are interested in understanding how the connections between different dynamical subsystems enable or jeopardize the structural controllability of the overall system.  This study follows the general lines provided in~\cite{StructNet}, but the verification procedures proposed were based in matrix nets leading to a computational effort that increases exponentially with the dimension of the problem.  Alternatively,  an efficient method is proposed in~\cite{Davison1977109} that accounts for the whole system instead of local properties (i.e., the components of the system and their interconnections), but this method does not apply to arbitrary systems.  More precisely, it is assumed that when connected, the state space digraph (see Section~\ref{problemstatement} for formal definition) is spanned by a disjoint union of cycles, which is called a \emph{rank constraint}.  In contrast, in~\cite{Rech1991877} and~\cite{gcactus}, the authors have presented results on the structural controllability of interconnected dynamical systems with interconnection in cascade.  Nevertheless, these structures are not unique, and the interconnection of these is established assuming such connectible structures are given; therefore, no practical criteria for computing such structures and verifying the results were provided.  More recently, in~\cite{blackhall} similar results were obtained by exploring which variables may belong to a structure and referred to as controllable state variables.  Thus, similarly to~\cite{Rech1991877} and~\cite{gcactus}, the results depend on the identified structures, but no method to systematically identify these structures is provided.  In~\cite{Yang19951011} the study is conducted assuming that all the subsystems, except a \emph{central} subsystem which is allowed to interact with every other subsystem, have the same dynamic structure and the interconnection between the several subsystems also has the same structure (even though they may not be used).  The implications of local interactions on the system's structural controllability determined in~\cite{Yang19951011}  can be obtained as particular solutions to the results  proposed hereafter.

In~\cite{PequitoJournal}, we studied the problem of determining the sparsest input matrix ensuring structural controllability, and additionally addressed the sparsest co-design problem, i.e., determining the sparsest solution of inputs, outputs and information patterns such that decentralized control laws guaranteeing that the closed-loop system enables arbitrary pole placement, based in static output feedback, exist.  Further, polynomial algorithms with computational complexity $\mathcal O(n^3)$ were provided to both problems, where $n$ is the the number of state variables.  In~\cite{NPcompStuct}, we extended the results in~\cite{PequitoJournal}, to the setting where the selection of inputs is constrained to a given collection of inputs, and sparsity of the inputs is not accounted for.  And later, in~\cite{PequitoCost} the problem in \cite{PequitoJournal} was further extended to determining the input matrix incurring in the minimum cost with respect to some cost function, and ensuring structural controllability.  Further, procedures with $\mathcal O(n^{\omega})$ computational complexity where provided, where $\omega<2.373$ is the lowest known exponent  associated with the complexity of multiplying two $n\times n$ matrices.  All these contrast with the one addressed in the current paper in the sense that we aim to verify  structural controllability properties in a distributed fashion.

On the other hand, composability aspects regarding controllability have been heavily studied by several authors, see for instance,~\cite{Zhou201563,ChenComposability,Wolovich1974209, Yonemura, Davison73, DavisonNewCompose}.  Briefly, all these resort to the well known Popov-Belevotch-Hautus (PBH) eigenvalue controllability criterion for LTI systems~\cite{Hespanha09}.  We notice that this criterion requires the knowledge of the overall system to infer its controllability.  The reason is closely related with the loss of degrees of freedom imposed by interconnected dynamical systems, as well as conservation laws in general, that result in the reduction of the rank of the system's dynamics matrix when compared with the sum of the ranks of the dynamics matrices of each subsystem.  Consequently, even if all subsystems are controllable, after their interconnection the resulting dynamical system  may not be.  Notwithstanding, the same does not happen when dealing with structural systems, where if all subsystems are structurally controllable, then the overall system is structurally controllable.  This is why, the key technical results in this paper correspond to verifying structural controllability in scenarios where the individual subsystems are not structurally controllable by themselves. In addition, while not guaranteeing that a system is controllable, we can regard the structural controllability analyses in this paper as necessary conditions for controllability.

\vspace{0.2cm}
\noindent \textbf{Main Contributions: } The main contributions of this paper are threefold:

\noindent $(i)$  We provide sufficient conditions for similar systems to be structurally controllable.  More precisely, these rely only on the structure of the subsystem and interconnection between subsystems.  A distributed algorithm is proposed, that can verify these conditions in polynomial time.

\noindent $(ii)$  We provide sufficient conditions for serial systems to be structurally controllable.  A distributed algorithm to verify these conditions is provided.  It requires only the knowledge of the subsystem and its neighbors' structure, as well as its interconnections.  This algorithm requires that a subsystem has  the capability  to interact with its neighbors only, and has computational complexity equal to $\cO\left(n_S^{3}\right)$, where $n_S$ denotes the total number of state variables and inputs present in a subsystem and its neighbors.

\noindent $(iii)$  We provide a distributed algorithm to verify necessary and sufficient conditions to ensure structural controllability for any interconnected dynamical system that consists of LTI subsystems.  This algorithm requires that a subsystem has the capability  to interact with its neighbors only, have access to its own structure and partial information regarding decisions performed by its neighbors that do not require sharing (directly) the structure of the neighboring agents.

The rest of this paper is organized as follows.  In Section~\ref{problemstatement} we formally describe the problem statement.  Section~\ref{prelim}  introduces some concepts in structural systems theory, that will be used throughout the remainder of the paper.  The main contributions are presented in Section~\ref{mainresults}  and in Section~\ref{illustrativeexample} we provide examples that illustrate the main findings.  Finally, Section~\ref{conclusions} concludes the paper and discusses avenues for further research.

\section{Problem Statement}
\label{problemstatement}

Consider $r$ linear time-invariant (LTI) dynamical systems described by
\[
    \dot x_i(t) =A_ix_i(t)+B_iu_i(t),\quad i=1,\ldots,r,
\]

\noindent where $x_i\in\mathbb{R}^{n_i}$ is the state, and $u_i\in\mathbb{R}^{p_i}$ the input of the $i$-th system.  The dynamical system can be described by the pair $(A_i,B_i)$, where $A_i\in\mathbb{R}^{n_i\times n_i}$ is the dynamic matrix of subsystem $i$ and $B_i\in\mathbb{R}^{n_i\times p_i}$ its input matrix.  By considering the interconnection from subsystem $i$ to subsystem $j$ for all possible subsystem  pairs we obtain the interconnected dynamical system described as follows:

{%
    \footnotesize
    \begin{equation}
        \dot x(t) =\underbrace{%
            \begin{bmatrix}
                A_1 &  E_{1,2} & \cdots &  E_{1,r} \\
                E_{2,1} &  A_2 & \cdots& E_{2,r} \\
                \vdots & \ddots& \ddots& \vdots \\
                E_{r,1} &\cdots &E_{r,r-1} & A_r
            \end{bmatrix}
        }_Ax(t) +
        \underbrace{%
            \begin{bmatrix}
                B_1 & 0 & \cdots & 0 \\
                0 & B_2 & \cdots & 0 \\
                \vdots & \ddots & \ddots & \vdots \\
                0 & 0 & \cdots & B_r
            \end{bmatrix}
        }_Bu(t),
        \label{eq:inputDynamics}
    \end{equation}
}

\noindent
where the state is given by $x=[x_1^\tr\ \dots \ x_r^\tr]^\tr\in \R^n$,  with $n=\sum_{i=1}^r n_i$, and the input given by $u=[u_1^\tr\ \dots \ u_r^\tr]^\tr\in \R^p$, with $p=\sum_{i=1}^r p_i$.  In addition,  $E_{i,j}\in\R^{n_i\times n_j}$ is referred to as the \emph{connection matrix} from the $j$--th subsystem to the $i$--th subsystem.  Further, we denote the system~\eqref{eq:inputDynamics} by the matrix pair $(A,B)$, denoting the $i$--th subsystem, $i=1,\dots,r$ of~\eqref{eq:inputDynamics} by the matrix pair $(A_i,B_i)$.  Finally, we call those subsystems $(A_j,B_j)$, with $j=1,\dots,r$ such that $E_{j,i}\neq 0$, the \emph{outgoing neighbors} of the $i$--th subsystem, and those that $E_{i,j}\neq 0$ the \emph{incoming neighbors} of the $i$--th subsystem; we refer to them collectively as the \emph{neighbors} of the $i$--th subsystem.

Now, consider the sparsity pattern of the matrix pair $(A,B)$ which we denote by the structural system $(\bar A,\bar B)$; similarly, we denote by $(\bar A_i,\bar B_i)$ the structural pair of matrices associated with $(A_i,B_i)$, and $\bar E_{i,j}$ the sparsity pattern of $E_{j,i}$.  Within the context of structural systems, the structural counterpart of controllability can be introduced as follows.

\begin{definition}[\cite{dionSurvey}]
    Given a structural system $(\bA,\bB)$, we say that it is structurally controllable if and only if, there exists at least one control system $(A,B)$ with the same sparsity pattern as $(\bA,\bB)$ (i.e.~$A_{i,j} = 0$ if $\bA_{i,j}=0$ and $B_{i,k}=0$ if $\bB_{l,k} =0$) which is controllable.
    \label{def:structural-control}
\end{definition}

It can be seen, from density arguments, that if $(\bA,\bB)$ is structurally controllable, then almost all control systems $(A,B)$ with the same sparsity as $(\bA,\bB)$ are structurally controllable~\cite{dionSurvey}.  We say that a control system $(A,B)$ is \emph{structurally controllable} if the associated structural system $(\bA,\bB)$ is structurally controllable.

Therefore the problem addressed in the current paper can be posed as follows.

\subsubsection*{Problem}
Given a collection of control systems $(A_i,B_i)$,  $i = 1,\dots,r$, and the interconnection from each subsystem $i$ to its neighbors, i.e.,  $(A_j,B_j,E_{j,i})$ for all $j\neq i$, design a distributed procedure to determine if the interconnected control system $(A,B)$ given in~\eqref{eq:inputDynamics} is structurally controllable.  \hfill $\circ$

Further, note that in a non-structural setting, local properties are not enough to guarantee controllability since the connection to other subsystems may lead to parameter cancellations~\cite{Davison73,DavisonNewCompose}; therefore, by resorting to structural system theoretic study, the approach presented hereafter allows us to obtain only necessary conditions for controllability.

\section{Preliminaries and Terminology}\label{prelim}

In this section, we review some preliminary concepts used to analyze the problem of structural controllability of interconnected dynamical systems.

\subsection{Structural Systems}

In order to perform structural analysis efficiently, it is customary to associate to~\eqref{eq:inputDynamics} a directed graph, or \emph{digraph} $\cD=(\cV,\cE)$, in which $\cV$ denotes the set of \textit{vertices} and $\cE$ the set of \textit{edges}, where $(v_j,v_i)$ represents an \emph{edge from the vertex $v_j$ to the vertex $v_i$}.  To this end, let $\bA\in\{0,1\}^{n\times n}$ and $\bB\in\{0,1\}^{n\times p}$ be the binary matrices that represent the sparsity patterns of $A$ and $B$  as in~\eqref{eq:inputDynamics}, respectively.  Denote by $\cX=\{x_1,\cdots,x_n\}$ and $\cU=\{u_1,\cdots,u_p\}$ the sets of \emph{state} and \emph{input vertices}, respectively.  And by $\cE_{\cX,\cX}=\{(x_i,x_j):\ \bA_{ji}\neq 0\}$, $\cE_{\cU,\cX}=\{(u_j,x_i):\ \bB_{ij}\neq 0\}$, the sets of edges between the vertex sets in subscript.  We may then introduce the \emph{state digraph} $\cD(\bA)=(\cX,\cE_{\cX,\cX})$ and the \emph{system digraph} $\cD(\bA,\bB)=(\cX\cup \cU,\cE_{\cX,\cX}\cup \cE_{\cU,\cX} )$.  Note that in the digraph $\cD(\bA,\bB)$, the input vertices representing the zero columns of $\bB$ correspond to isolated vertices.  As such, the number of \emph{effective} inputs, i.e., the inputs which actually exert control, is equal to the number of nonzero columns of $\bB$, or, in the digraph representation, the number of input vertices that are connected to at least one state vertex through an edge in $\cE_{\cU,\cX}$.

A \emph{directed path} from the vertex $v_1$ to $v_k$ is a sequence of edges $\{(v_1,v_2),(v_2,v_3),\hdots,(v_{k-1},v_k)\}$.  If all the vertices in a directed path are distinct, then the path is said to be an \emph{elementary path}.  A \emph{cycle} is an elementary path from $v_1$ to $v_k$, together with an edge from $v_k$ to $v_1$.

Given a digraph $\cD=(\cV,\cE)$, we say that $\cD'=(\cV',\cE')$ is a subgraph of $\cD$ if it is a digraph with $\cV'\subseteq \cV$ and $\cE'\subseteq\cE$, which we denote by $\cD'\subseteq\cD$.  Further, we say that $\cD'$ spans $\cD$ if $\cV'=\cV$.

We also require the following graph-theoretic notions~\cite{Cormen}: A digraph $\cD$ is strongly connected if there exists a directed path between any two vertices.  A \emph{strongly connected component} (SCC) is a subgraph $\cD_S=(\cV_S,\cE_S)$ of $\cD$ such that for every $u,v \in\cV_S$ there exist paths from $u$ to $v$ and from $v$ to $u$ and is maximal with this property (i.e., any subgraph of $\cD$ that strictly contains $\cD_S$ is not strongly connected).

\begin{definition}[\cite{PequitoJournal}]\label{linkedSCC}
    An SCC is said to be linked if it has at least one incoming or outgoing edge from another SCC\@.  In particular, an SCC is \textit{non-top linked} if it has no incoming edges to its vertices from the vertices of another SCC.
    \hfill $\diamond$
\end{definition}

Further, given a digraph $\cD = (\cV,\cE)$, we say that $\cD$ is a \emph{weakly connected digraph} if $(\cV,\cE\cup\cE^\tr)$ is strongly connected, where $\cE^\tr \equiv \{ (v',v) : (v,v')\in\cE\}$

\begin{figure}[htpb]
    \centering
\includegraphics[scale=0.35]{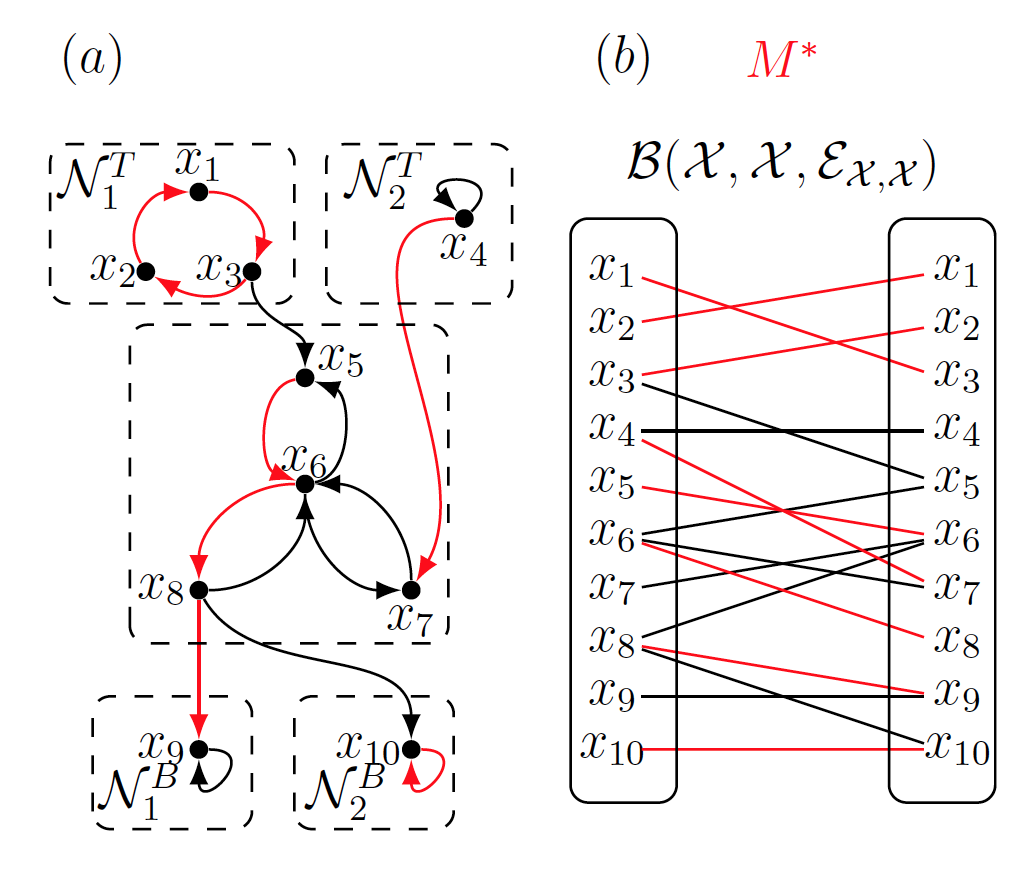}
    \caption{\footnotesize In $(a)$ we have depicted a digraph $\cD(\bA)$.  The SCCs are depicted within the dashed boxes; the digraph consists of two non-top linked SCCs ($\cN^T_1,\cN^T_2$) and two non-bottom linked SCCs ($\cN^B_1,\cN^B_2$).  In $(b)$ we depicted the state bipartite graph $\cB(\bA)\equiv \cB(\cX,\cX,\cE_{\cX,\cX})$, along with a \emph{maximum matching} $M^*$ comprising the edges in red, which are also highlighted in $(a)$ as they comprise a \emph{minimum path and cycle decomposition}, as explained in Lemma~\ref{lem:path-cycle-decomp}.}
    \label{fig:exampleBackground}
\end{figure}

For any digraph $\cD = (\cV,\cE)$ and any two sets $\cS_1, \cS_2\subset \cV$ we define the \textit{bipartite graph} $\cB(\cS_1,\cS_2,\cE_{\cS_1,\cS_2})$ 
where we call $\cS_1$ the set of \emph{left vertices}, and $\cS_2$ the set of \emph{right vertices}; and the edge set $\cE_{\cS_1,\cS_2}=\cE\cap(\cS_1\times\cS_2)$.  We call the bipartite graph $\cB(\cV,\cV,\cE)$ the bipartite graph associated with $\cD(\cV,\cE)$.  In the sequel we will make use of the \emph{state bipartite graph}, $\cB(\bA)\equiv \cB(\cX,\cX,\cE_{\cX,\cX})$, which is the bipartite graph associated associated with the state digraph $\cD(\bA)=(\cX,\cE_{\cX,\cX})$, and the \emph{system bipartite graph} $\cB(\bA,\bB) = \cB(\cU\cup\cX,\cX,\cE_{\cX,\cX}\cup\cE_{\cU,\cX})$.  An illustration of the state bipartite graph can be found in Figure~\ref{fig:exampleBackground}$-(b)$ associated with the state digraph in Figure~\ref{fig:exampleBackground}$-(a)$.

Given a bipartite graph $\cB(\cS_1,\cS_2,\cE_{\cS_1,\cS_2})$, a matching $M$ corresponds to a subset of edges in $\cE_{\cS_1,\cS_2}$ so that no two edges have a vertex in common, i.e., given edges $e=(s_1,s_2)$ and $e'=(s_1',s_2')$ with $s_1,s_1' \in \cS_1$ and $s_2,s_2'\in \cS_2$, $e, e' \in M$ only if $s_1\neq s_1'$ and $s_2\neq s_2'$.  An example of a matching is provided in Figure~\ref{fig:exampleBackground}$-(b)$.

A maximum matching $M^*$ is a matching $M$ that has the largest number of edges among all possible matchings, as the one shown in Figure~\ref{fig:exampleBackground}$-(b)$.

Further, it is possible to assign a weight to the edges in a bipartite graph, say $c(e)$ (where $c$ is a function from $\cE_{\cS_1,\cS_2}$ to $\R^+$).  We thus obtain a \emph{weighted bipartite graph}, and can introduce the concept of \emph{minimum weight maximum assignment} problem.  This problem consists in determining a maximum matching whose overall weight is as small as possible, i.e., a matching $M^c$ such that
\[
M^c=\arg\min_{M\in \cM} \sum_{e \in M} c(e),
\]
where $\cM$ is the set of all maximum matchings.  This problem can be efficiently solved using the Hungarian algorithm~\cite{Munkres1957},  with complexity of $\cO\left( \max\{|\cS_1|,|\cS_2|\}^3\right)$.  We call the vertices in $\cS_1$ and $\cS_2$ belonging to an edge in $M^*$, the \textit{matched vertices} with respect to (w.r.t.) $M^*$, otherwise, we call them \textit{unmatched vertices}.  It is worth noticing that there may exist more than one maximum matching, for example, in Figure~\ref{fig:exampleBackground}$-(b)$ replacing the edge $(x_8,x_9)$ in $M^*$ with $(x_9,x_9)$ yields a different matching with the same number of edges.
For ease of referencing, in the sequel, the term \emph{right-unmatched vertices}, with respect to $\cB(\cS_1,\cS_2,\cE_{\cS_1,\cS_2})$ and a matching $M$, not necessarily maximum, will refer to those vertices in $\cS_2$ that do not belong to an edge in $M^*$, and by duality a vertex from $\cS_1$ that does not belong to an edge in $M^*$ is called a \emph{left-unmatched vertex}.

The following result gives us some liberty in choosing the right- and left-unmatched vertices of the maximum matching of a bipartite graph.

\begin{lemma}[\cite{PequitoJournal}]\label{maxMatDecomp}
    Let $\cB(\cV,\cV,\cE)$ be a bipartite graph.  If $M^1$ and $M^2$ are two possible maximum matchings of $\cB(\cV,\cV,\cE)$ with right-unmatched and left-unmatched vertices given by $(\cU_R(M^1),\cU_L(M^1))$ and $(\cU_R(M^2),\cU_L(M^2))$ respectively, then, there exists a maximum matching $M^*$ of $\cB(\cV,\cV,\cE)$ with sets of right-unmatched and left-unmatched vertices given by $(\cU_R(M^1)$, $\cU_L(M^2))$.
    \hfill $\diamond$
\end{lemma}

Now, we can interpret a maximum matching of a bipartite graph associated to a digraph, at the level of the digraph.

\begin{lemma}[Maximum Matching Decomposition~\cite{PequitoJournal}]\label{lem:path-cycle-decomp}
    Consider the digraph $\cD=(\cV,\cE)$ and let $M^*$ be a maximum matching associated with the bipartite graph $\cB(\cV,\cV,\cE)$.  Then, the digraph $\cD=(\cV, M^*)$  comprises a disjoint union of cycles and elementary paths (by definition an isolated vertex is regarded as an elementary path with no edges), beginning in the right-unmatched vertices and ending in the left-unmatched vertices of $M^*$, that span $\cD$.  Moreover, such a decomposition is \emph{minimal}, in the sense that no other spanning subgraph decomposition of $\cD(\bA)$ into elementary paths and cycles contains strictly fewer elementary paths.
    \hfill $\diamond$
\end{lemma}

In addition, to make comparisons with previous work (namely, \cite{Rech1991877} and \cite{gcactus}), we need to introduce the following definitions.

\begin{definition}[\cite{Lin_1974}]\label{def:cactus}
    Given a digraph $\cD$, an elementary path in $\cD$, also called a \emph{stem}, is a cactus.  Given a cactus $\cG=(\cV_\cG,\cE_\cG)\subseteq\cD$, and a cycle $\cC=(\cV_\cC,\cE_\cC)\subseteq\cD$, such that $\cG$ and $\cC$ have no vertices in common, and there is an edge from a vertex in $\cG$ to a vertex in $\cC$, then $\cG\cup\cC = (\cV_\cG\cup\cV_\cC , \cE_\cG\cup\cE_\cC)$ is a cactus.
    \hfill$\diamond$
\end{definition}

Particularly, in the case where $\cD = \cD(\bA,\bB)$, a cactus $\cG$ in $\cD$ is called an \emph{input cactus} if the stem starts at an input vertex.  Further, we note that the decomposition into disjoint elementary paths and cycles in Lemma~\ref{lem:path-cycle-decomp}, can be used to determine a spanning decomposition of the graph into disjoint cacti~\cite{PequitoJournal}.

When dealing with interconnected dynamical systems, the structure of the connection between the subsystems will create connections between the SCCs of different subsystem digraphs.  This, in turn, makes it difficult to identify the SCCs of the system digraph of the overall system by analysing the SCCs of each subsystem digraph seperately and the connections to their neighbors.  Hence, we introduce the concept of \emph{reachability}~\cite{dionSurvey}.  We say that a state vertex $x$ in a system digraph is \emph{input-reachable} or \emph{input-reached} if there exists a path from an input vertex to it.

All of these constructions can be used to verify the structural controllability of an LTI system by analysing the associated graphs.

\begin{theorem}[\cite{dionSurvey,PequitoJournal}]\label{thm:struct-cont-cond}
    For LTI systems described by \eqref{eq:inputDynamics}, the following statements are equivalent:
    \begin{enumerate}
        \item[$(1)$] The corresponding structured linear system $(\bA,\bB)$ is structurally controllable;
        \item[$(2)$] The digraph $\cD(\bA,\bB)$ is spanned by a disjoint union of input cacti;
        \item[$(3i)$]\label{nontopsccinp}  The non-top linked SCCs of the system digraph $\cD(\bA,\bB)$ are comprised of input vertices, and
        \item[$(3ii)$]\label{norightunmatch} there is a matching of the system bipartite graph $\cB(\bA,\bB)$ without right-unmatched vertices;
        \item[$(4i)$]\label{reachedinp} Every state vertex is input-reachable, and
        \item[$(4ii)$]\label{norightunmatch2} there is a matching of the system bipartite graph $\cB(\bA,\bB)$ without right-unmatched vertices.
    \end{enumerate}
            \hfill$\diamond$
\end{theorem}

\section{Main Results}
\label{mainresults}

\begin{figure}[htb]
    \centering
\includegraphics[scale=0.28]{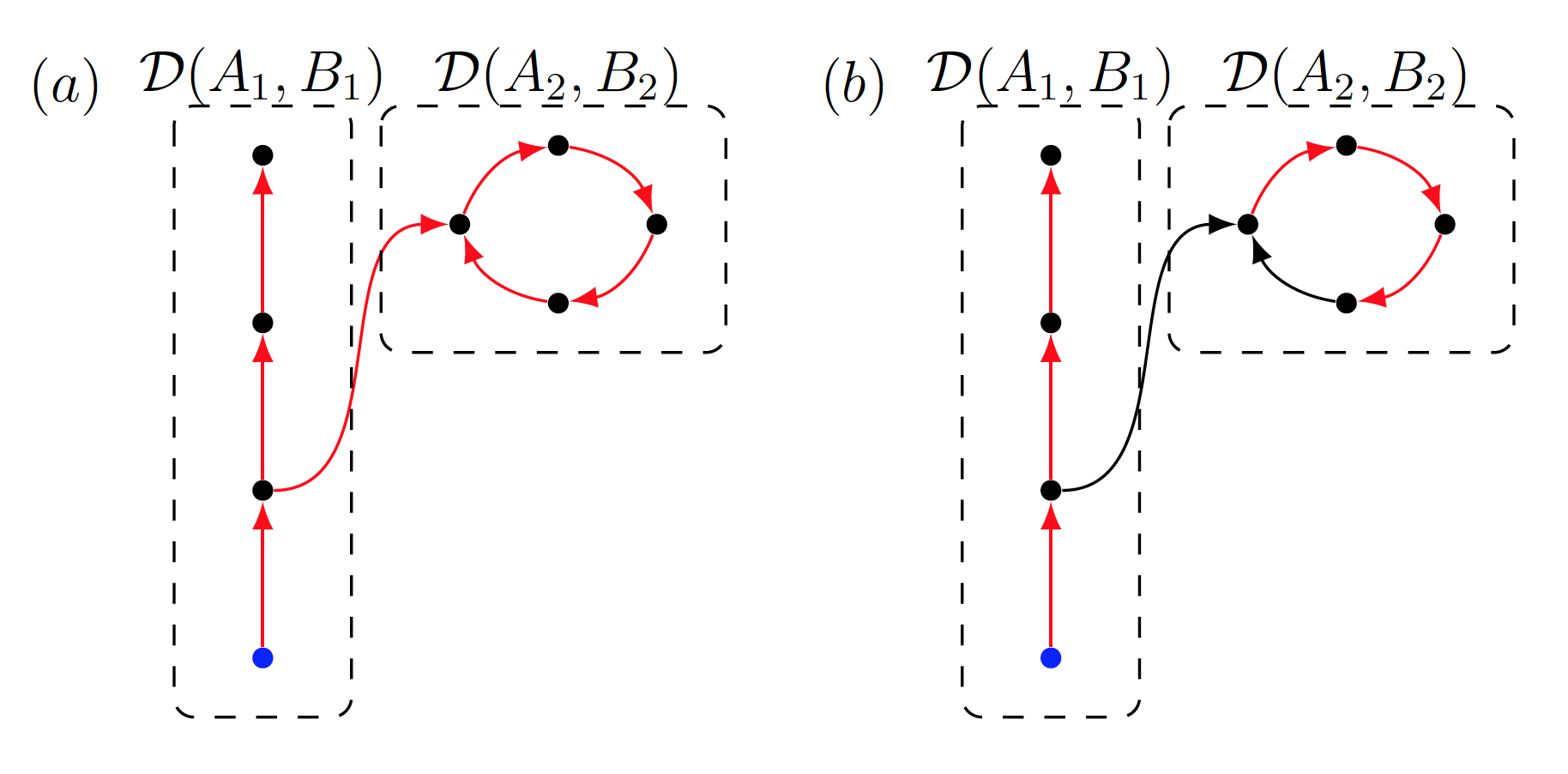}
    \caption[]{\footnotesize In $(a)$, we present the digraph associated to an interconnected dynamical system, where the different subsystems are represented inside the dashed boxes.  Recall the definition of cactus, it can be readily be seen that the digraph $\cD(\bA,\bB)$ is spanned by the input cactus, depicted by the red edges, rendering the structural system $(\bA,\bB)$ structurally controllable by Theorem~\ref{thm:struct-cont-cond}--$(2)$.  In $(b)$, however, we depict possible cacti that span each of the subsystem digraphs.  Since a spanning cactus for $\cD (\bA_2,\bB_2)$ has to include a stem comprising at least a vertex, neither of the cacti spanning $\cD(\bA_1,\bB_1)$ and $\cD(\bA_2,\bB_2)$ can contain any cycles, and there is no way of prolonging the stem that spans $\cD(\bA_1,\bB_1)$ to include a stem spanning $\cD(\bA_2,\bB_2)$.  This shows that the conditions proposed in~\cite{Rech1991877} are not necessary.}
    \label{fig:nocactus}
\end{figure}

We begin this section by providing sufficient conditions for an interconnected dynamical system to be structurally controllable in the case where all the subsystems have the same structure (Theorem~\ref{thm:perf-match-sim-system} and Theorem~\ref{thm:right-unmatch-sim-syst}).  We then focus on more general interconnected dynamical systems, called \emph{serial systems}, and propose sufficient conditions for their structural controllability (Lemma~\ref{lem:sufP1}); in addition, an efficient distributed algorithm (Algorithm~\ref{alg:sufP1}) to verify these conditions is provided, which has its correctness and complexity proven in Theorem~\ref{thm:sufP1alg}.  In light of these conditions, we explain why previous results in this line~\cite{Rech1991877} presented conditions that are only sufficient instead of necessary and sufficient (Figure~\ref{fig:nocactus}).  Finally, we end this section by describing an efficient distributed algorithm (Algorithm~\ref{alg:fdic}) to verify the structural controllability of an arbitrary interconnected dynamical system, which has its correctness and complexity proven in Theorem~\ref{thm:algc-fdic}. In order to perform this verification, each subsystem has to perform calculations using the information about itself and its neighbors.

Often, interconnected dynamical systems under analysis are composed of subsystems that are similar, as formally introduced next.

\begin{definition}
    Let $\bE\in\{0,1\}^{r\times r}$,  $\bA',\bH\in\{0,1\}^{n\times n}$, \mbox{ $\bB'\in\{0,1\}^{n\times p}$}, be matrices with the restriction that $\bE_{i,i}=0$ ($i=1,\dots,r$).  Then, we denote by $(\bA',\bB',\bH,\bE')$ the structural system $(\bA,\bB)$, where $\bA = (I_r\otimes\bA')\lor(\bE\otimes\bH)$ and $\bB = I_r\otimes\bB'$, where $\otimes$ denotes the matrix Kronecker product.  Further, we say this system to be \emph{composed of $r$ similar subsystems}, or, in short, a \emph{similar system}.
    \hfill $\diamond$
    \label{def:similar-systems}
\end{definition}

\begin{remark}
    Note that in the case of similar systems, $\bH$ is the structural matrix modeling the interactions between each subsystem and its neighbors, all of which have the same structure.
\end{remark}

\begin{definition}
    Let $(\bar A,\bar B)$ be the structural matrices associated with the interconnected dynamical system in~\eqref{eq:inputDynamics}. We define the \emph{condensed graph of the system} as being the digraph $\cD^*(\bA)\equiv\cD(\cA,\cE)$, where $a_i \in \cA \equiv \{a_1,\dots,a_r\}$ is a vertex representing the $i$--th subsystem, and $(a_i,a_j)\in \cE \equiv \{(a_i,a_j) | E_{j,i}\neq 0\}$ a directed edge representing a communication from subsystem $j$ to subsystem $i$.  Moreover, if there is no directed edge ending in a vertex, this vertex is referred to as a \emph{source}.
    \hfill $\diamond$
    \label{def:condensed-graph}
\end{definition}

Note that in the case that a system $(\bA,\bB)$ is composed of $r$ similar systems and parametrized by matrices $(\bA',\bB',\bH,\bE)$ the condensed graph $\cD^*(\bA)$ is the same as the digraph $\cD(\bE)$.  Now, we assess the structural controllability of these systems when their subsystems are not structurally controllable.

\begin{theorem}
    Let the system $(\bA,\bB)$ be composed of $r$ similar components, and parametrized by $(\bA',\bB',\bH',\bE)$, where $(\bA',\bB')$ is not structurally controllable, and $\cB(\bA',\bB)$ has a matching without right-unmatched vertices.  The pair $(\bA,\bB)$ is structurally controllable if and only if $(\bA'\lor\bH,\bB')$ is structurally controllable and the condensed graph $\cD^*(\bA)$ has no sources.
    \hfill$\diamond$
    \label{thm:perf-match-sim-system}
\end{theorem}

\begin{IEEEproof}
    To prove the equivalence, we begin by proving that the conditions are sufficient by contrapositive; subsequently, we prove directly that the conditions are also necessary.

    First, notice that since $(\bA',\bB')$ is not structurally controllable despite $\cB(\bA',\bB')$ having a maximum matching without right-unmatched vertices, it follows, from Theorem~\ref{thm:struct-cont-cond}--$(4)$, that $\cD(\bA',\bB')$ has a vertex which is not reachable from any input vertex.  Subsequently, assume that $\cD^*(\bA)$ has a source, then it follows that there is a subsystem $(\bA'_j,\bB'_j)$ with no incoming edges from other subsystems, and so the overall system digraph $\cD(\bA,\bB)$ has a state vertex without a path from any input vertex to it. Hence, by Theorem~\ref{thm:struct-cont-cond}--$(4)$, $(\bA,\bB)$ is not structurally controllable.  Further, $(\bA'\lor\bH,\bB')$ is not structurally controllable.  Since $\cB(\bA',\bB')$ has a matching without right-unmatched vertices, so does $\cB(\bA'\lor\bH,\bB')$, which means, by Theorem~\ref{thm:struct-cont-cond}--$(4)$ that, $\cD(\bA'\lor\bH,\bB')$ must have a state vertex which is not reachable from any input vertex.  This implies that the corresponding state vertex is not reachable from an input vertex in any of the subsystems (since a path from an input vertex in the overall system translates into one such path in $\cD(\bA'\lor\bH,\bB')$).

    Finally, assume that $(\bA'\lor\bH,\bB')$ is structurally controllable and that $\cD^*(\bA)$ has no sources. Then, for each state vertex, there is a path from an input vertex to it, which implies, by Theorem~\ref{thm:struct-cont-cond}--$(4)$, that $(\bA,\bB)$ is structurally controllable.
\end{IEEEproof}

In the next result, we relax the assumptions from Theorem~\ref{thm:perf-match-sim-system}, about the structure of the dynamics of the subsystems.  Thus, allowing for applications in the design of interconnections between subsystems that may fail to meet such criteria.

\begin{theorem}
    Given an interconnected dynamical system $(\bA,\bB)$ composed of $r$ similar components, and parametrized by $(\bA',\bB',\bH,\bE)$, where $(\bA',\bB')$ is not structurally controllable, then $(\bA,\bB)$ is structurally controllable if both $(\bA'\lor\bH,\bB')$ is structurally controllable and $\cD^*(\bA)$ is spanned by cycles.
    \hfill$\diamond$
    \label{thm:right-unmatch-sim-syst}
\end{theorem}

\begin{IEEEproof}
    First, notice that if the digraph $\cD^* (\bA)$ is spanned by cycles, every vertex belongs to a cycle, and, in particular, it means that $\cD^* (\bA)$ has no sources.  In this case, the method of proof of Theorem~\ref{thm:perf-match-sim-system} is applicable to show that every state vertex has a path from an input vertex to it, so all that remains to show is that the $\cB(\bA,\bB)$ has no right-unmatched state vertices with respect to some maximum matching.  To this end, we first assume (without loss of generality) that $\cD^*(\bA)$ has one spanning cycle, and that the subsystems $(\bA'_1,\bB'_1),\cdots,(\bA'_r,\bB'_r)$ are ordered in such a way that $\bE_{i+1,i} = 1$ for $i=1,\dots,r-1$, and $\bE_{1,r}= 1$.

    Now, denote the state and input vertices of the $i$--th subsystem by $x_k^i$ with $k=1,\cdots,n$ and $u_l^i$ with $l=1,\cdots,m$, respectively.  In addition, let $M$ be a maximum matching of $\cB(\bA'\lor\bH,\bB')$ without right-unmatched state vertices, then we can partition $M'$ into three matchings $M'_B,M'_A,M'_H$ comprising, respectively.  The edges in $M'$ are of the form $(u_l,x_k)$,  $(x_l,x_k)$ when $\bA_{k,l}=1$, and the remaining ones are of the form $(x_l,x_k)$, when $\bA_{k,l}=0$ and $H_{k,l}=1$.  Finally, consider the matching $M$ of $\cB(\bA,\bB)$ comprising the edges:
    \begin{itemize}
        \item $(u^i_k,x^i_l)$, if $(u_k,x_l)$ is in $M'_B$;
        \item $(x^i_k,x^i_l)$, if $(x_k,x_l)$ is in $M'_A$;
        \item $(x^r_k,x^1_l)$, if $(x_j,x_l)$ is in $M'_H$;
        \item $(x^i_k,x^{i+1}_l)$, if $(x_j,x_l)$ is in $M'_H$.
    \end{itemize}
    To show that $M$ has no right-unmatched vertices associated with it, consider a state vertex $x^i_k$ of $\cB(\bA,\bB)$, and since $M'$ has no right-unmatched vertices, $x_k$ is not right-unmatched in $\cB(\bA',\bB')$; thus, there is an edge $(x_l,x_k)$ for some $l$ or $(u_{l'},x_k)$ in $M'$, but this implies that either $(x^i_l,x^i_k)$, $(x^{i-1}_l,x^i_k)$, or $(u^i_{l'},x^i_k) \in M$ (where $i-1 = r$ when $i=1$) from the construction above, thus $\cB(\bA,\bB)$ has no right-unmatched vertices w.r.t. $M$.

    Finally, if more than one cycle is necessary to span the digraph $\cD^*(\bA)$, then the same argument applies to each cycle individually.
\end{IEEEproof}

\begin{figure}[htb]
    \centering
\includegraphics[scale=0.25]{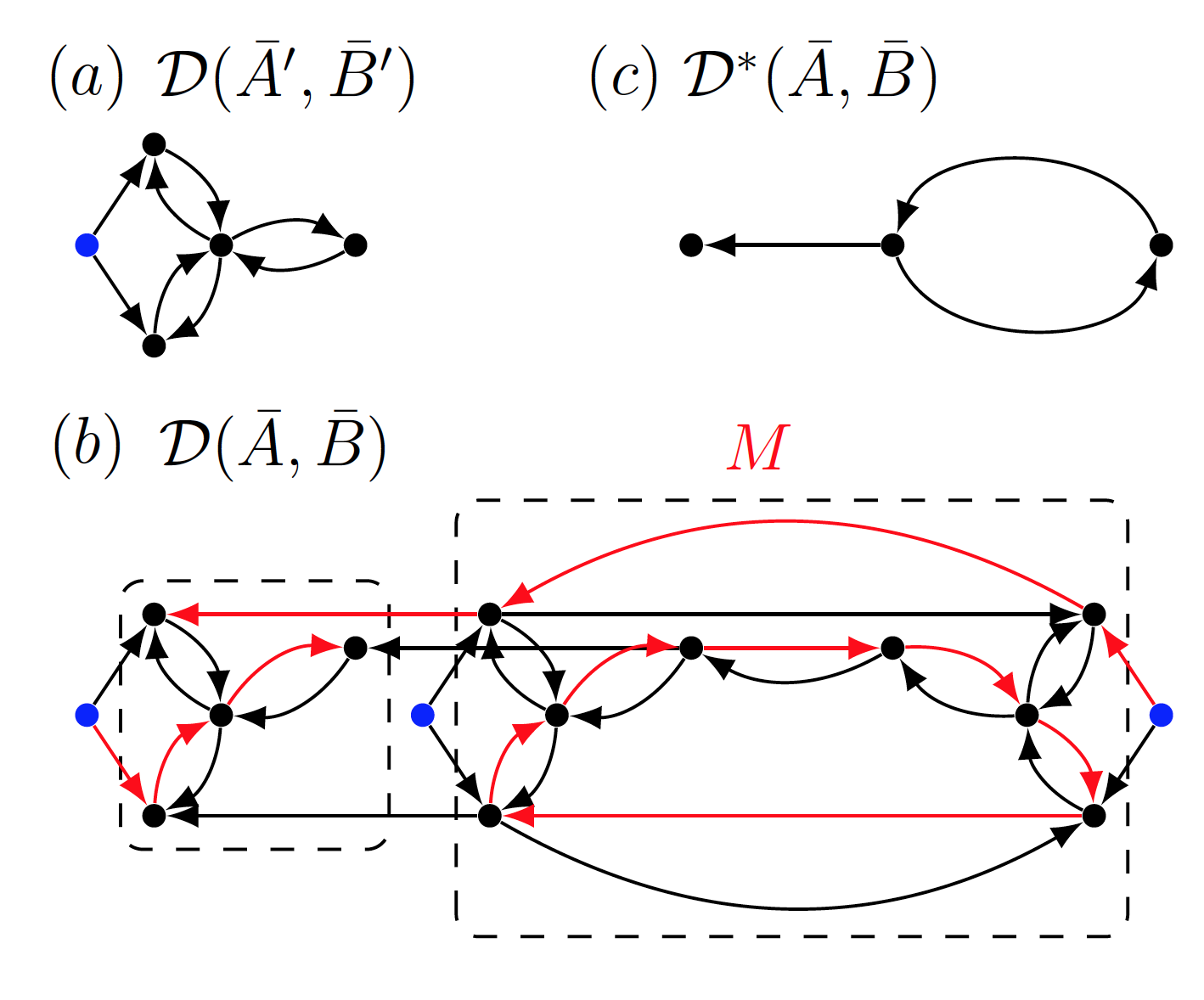}
    \caption{\footnotesize In $(a)$ we provide the digraph $\cD(\bA',\bB')$ of a system that is not structurally controllable, where we represent in blue the single input vertex.  By connecting three of these systems together as in $(b)$ the system $(\bA,\bB)$ becomes structurally controllable, as evidenced by the fact the matching $M$, associated to the path and cycle decomposition (Lemma~\ref{maxMatDecomp}) depicted by the red edges, has no right-unmatched state vertices (since every state vertex has an incoming red edge), and the fact that every non-top linked SCC is comprised of an input vertex.  Finally, in $(c)$  we show that the condensed graph $\cD^*(\bA)$ of the system in $(b)$ is not spanned by cycles, showing that the condition in Theorem~\ref{thm:right-unmatch-sim-syst} is not necessary.}
    \label{fig:example-not-necessary}
\end{figure}

\begin{remark}
    The conditions in Theorem~\ref{thm:right-unmatch-sim-syst} are not necessary, consider the example in Figure~\ref{fig:example-not-necessary}$-(b)$, that is structurally controllable, yet the condensed graph is not spanned by cycles.
    
    \hfill $\diamond$
\end{remark}

Now, we move toward methods for verifying structural controllability of interconnected dynamical systems by resorting to distributed algorithms.  To this end, we begin by introducing a result that will allow us to infer structural controllability of a family of interconnected dynamical systems that we refer to as serial systems.  Later, we provide a computational method to perform this verification in a distributed manner, that is, in which each subsystem only needs to have partial information about the system in order to verify if the system is structurally controllable or not.

It is worth noting that in order to ensure that all of the algorithms in the present paper can work as intended, several assumptions have to be made about the subsystems, and how they interact.  Namely, that each subsystem has a processing unit, and can send arbitrary messages to its neighboring subsystems; in addition, each subsystem is aware of the number of subsystems in the overall system and possesses a unique \emph{id}, and that the condensed graph of the system $\cD^*(\bA)$ is weakly connected.

\begin{lemma}
    Consider the structural system $(\bA,\bB)$ as in~\eqref{eq:inputDynamics} with subsystems $(\bA_1,\bB_1),\dots,(\bA_r,\bB_r)$.
    Then the system $(\bA,\bB)$ is structurally controllable, if there exist maximum matchings $M_0,\dots,M_r$ of the bipartite graphs $\cB(\bA_1),\dots,\cB(\bA_r)$ respectively, such that the following conditions hold:
    \begin{enumerate}
        \item \label{it:scc-top-linking} For each subsystem $(\bA_j,\bB_j)$, $j = 1,\dots,r$, the non-top linked SCCs of $\cD(\bA_j,\bB_j)$ consist of input vertices; and,
        \item \label{it:match-fitting} the bipartite graph 
            $$\cB\left(\bigcup_{i=1}^r\cU_L(M_i),\bigcup_{i=1}^r\cU_R(M_i),\bigcup_{i=1}^r\bigcup_{j\neq i}\cE_{\cU_L(M_j),\cU_R(M_i)}\right),$$
            admits a perfect matching, where $\cU_L(M_i)$ and $\cU_R(M_i)$ are the sets of left- and right-unmatched vertices of $M_{i}$ respectively, and $\cE_{\cU_L(M_j),\cU_R(M_i)}\subseteq\cE_{\cX,\cX}$ is the set of edges from vertices in $\cU_L(M_j)$ to vertices in $\cU_R(M_i)$.  \hfill $\diamond$
    \end{enumerate}
    \label{lem:sufP1}
\end{lemma}

\begin{IEEEproof}
    First, note that the non-top linked SCCs of $\cD(\bA,\bB)$ consist of SCCs of the subsystem digraphs $\cD(\bA_i,\bB_i)$. Further, for one such SCC to be non-top linked, it must contain at least one non-top linked SCC of some $\cD(\bA_i,\bB_i)$ in it.  Therefore, since every non-top linked SCC of every $\cD(\bA_i,\bB_i)$ is comprised of input vertices, and there are no edges from any neighboring system to input vertices, the non-top linked SCCs of $\cD(\bA,\bB)$ must contain input vertices.

    Secondly, note that the union of the maximum matchings $M_i$ of the $\cB(\bA_i,\bB_i)$ comprises a matching $M$ of $\cB(\bA,\bB)$.  Further, let $M'$ be the matching mentioned in the second condition.  Then, $M'$ is comprised of edges from left-unmatched vertices to right-unmatched vertices of $M$, and $M\cup M'$ is a matching of $\cB(\bA,\bB)$; further, since by hypothesis the matching $M'$ has no right-unmatched vertices, neither does $M\cup M'$.  Consequently, by Theorem~\ref{thm:struct-cont-cond}--$(3)$, this implies that the system is structurally controllable.
\end{IEEEproof}

Note that  using Lemma~\ref{lem:sufP1} we conclude that the system depicted in Figure~\ref{fig:nocactus}--$(a)$ is structurally controllable, yet by using the characterization in~\cite{Rech1991877}, it is not possible to obtain the same conclusion.  Further, Lemma~\ref{lem:sufP1} provides only a sufficient condition for structural controllability.  Nonetheless, these conditions can be verified in a distributed manner in the class of interconnected dynamical systems formally introduced next.

\begin{definition}
    We say that an interconnected dynamical system $(\bA,\bB)$ as in~\eqref{eq:inputDynamics} is a \emph{serial system} if each vertex of the condensed graph $\cD^*(\bA)$ has at most one outgoing edge.
    \label{def:serial-systems}
    \hfill $\diamond$
\end{definition}
\vspace{0.2cm}

\begin{figure}[htb]
    \centering
\includegraphics[scale=0.3]{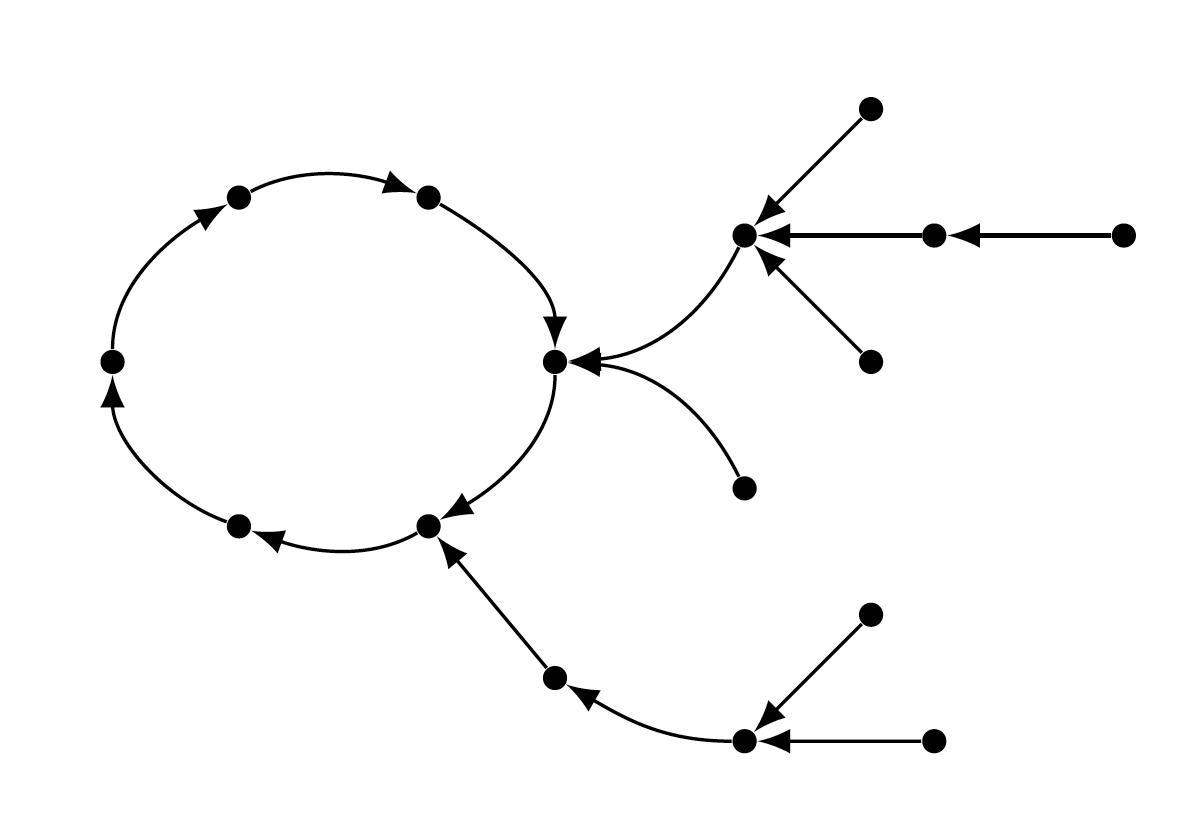}
    \caption{\footnotesize Example of a possible condensed graph for a serial system, where each vertex represents a subsystem, and each directed edge a non-zero connection matrix, see Definition~\ref{def:serial-systems}.}
    \label{fig:serial-system}
\end{figure}

Although serial systems seem a restrictive class of systems, they may exhibit a rich structure as  the condensed graph of a serial system   in Figure~\ref{fig:serial-system} illustrates.  Further, recall that serial systems enable us to verify the sufficient conditions for structural controllability in Lemma~\ref{lem:sufP1} in a distributed manner.  Thus, in Algorithm~\ref{alg:sufP1}, we present the procedure that each agent employs in order to verify the conditions in Lemma~\ref{lem:sufP1}.

Before introducing Algorithm~\ref{alg:sufP1}, we explain the functions that are used throughout the algorithm.  All these functions should be able to be applied by the $i$--th subsystem, $\textsc{Send}(x,j)$ sends the value of the variable $x$ to the $j$--th subsystem, while $\textsc{Rcv}(j)$ makes the system wait to receive a message from the $j$--th subsystem and subsequently reads this message.  Note that both these functions can only be applied when systems $i$ and $j$ interact with each other.  For communication to be successful, we assume that the systems perform these steps synchronously (i.e.\ they wait for the responses of their neighbors).  Finally, the procedure $\textsc{MinWtMaxMatch}(c,G)$ calculates the minimum weight maximum matching on the graph $G$ using  the cost function $c$, and boolean expressions contained in square brackets get evaluated (to $\textbf{True}$ or $\textbf{False}$).

\begin{algorithm}[!htpb]
    \caption{\footnotesize Distributed algorithm to verify sufficient conditions given in Lemma~\ref{lem:sufP1}, for an arbitrary serial system.}
    \label{alg:sufP1}
    \footnotesize
    \begin{algorithmic}[1]
        \Procedure{SeqStrtCtl}{$\bA_i,\bB_i,r,\bE_{i,j}\neq 0$}
        \LComment{$r = $ total number of subsystems, $\bE_{i,j} = $ connection matrices}
        \State $\texttt{nghI}(i) \gets \{ j : \bE_{i,j} \neq 0 \}$
        \State $\texttt{nghO}(i) \gets \{ j : \bE_{j,i} \neq 0 \}$
        \State $\texttt{nghbs}(i) \gets \texttt{nghI}(i) \cup \texttt{nghO}(i)$
        \LComment{send the dynamic matrix to (the unique) outgoing neighbor}
        \ForAll{$j\in\texttt{nghO}(i)$}
        \State $\Call{Send}{\bA_i,j}$
        \EndFor
        \LComment{receive the dynamic matrix of the incoming neighbors $j\in\texttt{nghI}(i)$}
        \ForAll{$j\in\texttt{nghI}(i)$}
        \Comment $\texttt{nghI}(i) = \{j_1,\dots,j_l\}$
        \State $\bA_j \gets \Call{Rcv}{j}$
        \EndFor
        \State $\bA'_i \gets \begin{bmatrix} \bA_i & \bE_{i,j_1} & \dots & \bE_{i,j_l}  \\ 0 & \bA_{j_1} & \dots & 0 \\ \vdots & \vdots & \ddots & \vdots \\ 0 & 0 & \dots & \bA_{j_l} \end{bmatrix}$
        \State $\bB_i'\gets \left[\bB_i^\tr , 0 , \dots, 0 \right]^\tr$
        \State $\cB_i \gets \cB(\bA'_i,\bB'_i)$
        \Function {$c_i$}{$y,x_k,$}
        \label{alg:sufP1:weight-func}
        \Comment define the weight function
        \If{$y=x_j, \text{ and } k,j \le n_i \text{ or } k,j > n_i$}
        \Statenl \Return $1$
        \Else
        \Statenl \Return $2$
        \EndIf
        \EndFunction
        \State $M_i \gets \Call{MinWtMaxMatch}{c_i,\cB_i}$
        \label{alg:sufP1:max-match}
        \State $\cU_R(M_i) \gets \{x_j:x_j\text{ right-unmatched w.r.t.} M_i\text{ and } j\le n_i\}$
        \label{alg:sufP1:right-unmatch}
        \State $\cN \gets \{x_j : x_j \text{ state vertex in non-top linked SCC of } \cD(\bA_i,\bB_i)\}$
        \label{alg:sufP1:unreached}
        \State $\texttt{mchd}(i) \gets [\cU_R(M_i) == \emptyset]$
        \label{alg:sufP1:mchd-ver}
        \State $\texttt{rchd}(i) \gets [\cN == \emptyset]$
        \label{alg:sufP1:rchd-ver}
        \LComment{check if the whole system can be structurally controllable according to whether the conditions are satisfied in the current system or not}
        \State $\texttt{ctld}(i) \gets \texttt{rchd}(i)\wedge\texttt{mchd}(i)$
        \label{alg:sufP1:guess}
        \For{$ k = 1 ,\dots, r$}
        \ForAll{$j\in\texttt{nghbs}(i)$}
        \label{alg:sufP1:consens-begin}
        \Statenl $\Call{Send}{\texttt{ctld}(i),j}$
        \Statenl $\texttt{ctld}(j) \gets \Call{Rcv}{j}$
        \EndFor
        \LComment{reconsider the answer in light of the values from the neighbors current answer}
        \State $\texttt{ctld}(i) \gets \texttt{ctld}(i)\bigwedge\limits_{j\in\texttt{nghbs}(i)} \texttt{ctld}(j)$
        \EndFor
        \label{alg:sufP1:consens-end}
        \LComment{return \textbf{True} if the system is structurally controllable and \textbf{False} otherwise}
        \State \Return $\texttt{ctld}(i)$
        \label{alg:sufP1:return}
        \EndProcedure
    \end{algorithmic}
\end{algorithm}

\begin{remark}
    The Algorithm~\ref{alg:sufP1} can be easily adapted to cover the case where each subsystem only has one incoming neighbor.  In this case, the only difference is that instead of having $\cB_i$ as the bipartite graph associated with the $i$--th subsystem and all its incoming neighbors, it would be that of the $i$--th subsystem and all its outgoing neighbors.
    \hfill$\diamond$
\end{remark}

The next result concerns the correctness and complexity of Algorithm~\ref{alg:sufP1}.

\begin{theorem}
    Algorithm~\ref{alg:sufP1} is correct, i.e.,~it verifies the sufficient conditions given in Lemma~\ref{lem:sufP1} for  serial systems.  Moreover, it has computational complexity $\cO\left(\max\limits_{i=1,\dots,r} N_i^3\right)$, with $N_i=m_i + \sum_{j\in \cI_i\cup\{i\}}n_j$, where $m_i$ and $n_i$ are the dimensions of the input and state space for the $i$--th subsystem, and $\cI_i\subseteq \left\{ 1,\dots,r \right\}$ is the set of subsystems which output is provided to the $i$--th subsystem.
    \hfill $\diamond$
    \label{thm:sufP1alg}
\end{theorem}

\begin{IEEEproof}
    To prove the correctness of Algorithm~\ref{alg:sufP1}, we start by proving the claim that a minimum weight maximum matching $M$ of $\cB_i$ w.r.t.\ the weight-function $c_i$ (defined in step~\ref{alg:sufP1:weight-func}) induces maximum matchings on $\cB(\bA_i)$, as well as on $\cB(\bA_j)$ for any subsystem $(\bA_j,\bB_j)$ with nonzero connection matrix $\bE_{i,j}$. Let $M_i$ be the matching resulting from restricting $M$ to the edges of $\cB(\bA_i)$, and in order to derive a contradiction, assume that $M_i$ is not a maximum matching of $\cB(\bA_i)$.  As a direct consequence of Berge's theorem (see for example Theorem~$1$ of~\cite{berge_two_1957}) the set of right-unmatched vertices of any matching contains the right-unmatched vertices of some maximum matching, so let $M'_i$ be a maximum matching such that $\cU_R(M'_i)\subseteq \cU_R(M_i)$.  Further, let $S_i\subseteq M$ be the set of edges from vertices not in $\cX_i$ to vertices in $\cX_i$, and let $S_i'$ be those edges in $S_i$ that end in vertices from $\cU_R(M'_i)$.  Now, $(M\setminus (M_i\cup S_i))\cup M'_i\cup S'_i$ is a matching of $\cB_i$ with the same number of edges as $M$ (since it has the same number of right-unmatched vertices) and with an overall weight lower than that of $M$ (since by hypothesis $S_i'\subsetneq S_i$), which contradicts the fact that $M$ is a minimum weight maximum matching.  The same argument works for the matching $M_j$ of $\cB(\bA_j)$ with $j\neq i$, replacing left-unmatched vertices with right-unmatched vertices.

    Now, since $(\bA,\bB)$ is a serial system, there is at most one $k\neq i$ with nonzero matrix $E_{k,i}$.  Therefore, we let $M'$ and $M''$ be the maximum matchings of $\cB(\bA_i)$ resulting from the maximum matchings of $\cB_i$ and $\cB_k$, respectively.  Then, by Lemma~\ref{maxMatDecomp}, there exists a maximum matching $M$ that has as left-unmatched vertices those of $M''$ and as right-unmatched vertices those of $M'$.  Subsequently, we only need to check for each subsystem that there is a minimum weight maximum matching of $\cB_i$ (w.r.t.~the weight function $w_i$) that has no right-unmatched state vertices; in summary, in Algorithm~\ref{alg:sufP1}, we set up the necessary structures until step~\ref{alg:sufP1:weight-func}.

    Now, in step~\ref{alg:sufP1:max-match}  the $i$--th subsystem computes the maximum matching $M_i$ of $\cB_i$, and in step~\ref{alg:sufP1:right-unmatch} the system calculates the associated set of right-unmatched vertices.  Next, in step~\ref{alg:sufP1:unreached} the subsystem calculates the set of state vertices in a non-top linked SCC of $\cD(\bA_i,\bB_i)$, and in steps~\ref{alg:sufP1:mchd-ver} and~\ref{alg:sufP1:rchd-ver}, it verifies the existence of right-unmatched state vertices of the $i$--th subsystem w.r.t.~the matching $M_i$, and the existence of in a non-top linked SCC of $\cD(\bA_i,\bB_i)$.  Finally, the subsystem decides if the whole system is structurally controllable or not in steps~\ref{alg:sufP1:guess}--\ref{alg:sufP1:consens-end}.  More precisely, after an initial conditions has been chosen and stored in $\texttt{ctld}(i)$, the subsystem updates this variable with the corresponding variable of its neighbors, and repeats this $r$ times.  Note that after $k$ iterations of the steps~\ref{alg:sufP1:consens-begin}--\ref{alg:sufP1:consens-end} the subsystem has updated $\texttt{ctld}(i)$ with the corresponding values of all subsystems at $k$ edges of distance from it.  Since the condensed graph of the systems is weakly connected, and there are only $r$ subsystems, $\texttt{ctld}(i)=\textbf{True}$ if and only if all subsystems had initially $\texttt{ctld}(j)=\textbf{True}$.  Finally, in step~\ref{alg:sufP1:return} the subsystem returns the value $\textbf{True}$ or $\textbf{False}$ depending on whether or not the system satisfies the conditions of Lemma~\ref{lem:sufP1}.

    Lastly, the complexity of Algorithm~\ref{alg:sufP1} is given as follows: since all of the steps have linear complexity except determining the minimum weight maximum matching of $\cB_i$ in step~\ref{alg:sufP1:max-match}, for which the Hungarian algorithm can be used with complexity $\cO\left(|N_i|^3\right)$, with $N_i = p_i + \sum_{j\in \cI_i\cup\{i\}}n_j$, $\bA_i \in \{0,1\}^{n_i\times n_i}$ and $\bB_i \in \{0,1\}^{n_i\times p_i}$ and $\cI_i\subseteq \left\{ 1,\dots,r \right\}$ is the set of indexes of subsystems incoming to the $i$--th subsystem~\cite{Munkres1957}.  This procedure has to be applied to each of the $r$ subsystems, which implies that the complexity of the algorithm becomes $\cO\left(\max\limits_{i=1,\dots,r} N_i^3\right)$.
\end{IEEEproof}

\begin{remark}
    If the system was not serial then there could be a subsystem, $k$ which outputs to both the $i$-- and $j$--th subsystems.  This  means that while computing maximum matchings of $\cB_i$ and $\cB_j$ separately, we can match the state vertex of the $k$--th subsystem to two different state vertices, one of the $i$--th subsystem and one of the $j$--th subsystem.  Further,  if there is a subsystem with incoming edges from every other system, the algorithm will calculate a maximum matching in a centralized manner.
\end{remark}

Now, we move toward distributed algorithms that can verify structural controllability of interconnected dynamical systems at large.  In this, each subsystem is required to share only partial information about its structure with its neighbors.  This algorithm, however, has a higher computational complexity than Algorithm~\ref{alg:sufP1}.  In order to infer structural controllability, we apply Theorem~\ref{thm:struct-cont-cond}--$(4)$, and begin by presenting an algorithm to verify if each of the state vertices, in the digraph associated to an interconnected dynamical system as in~\eqref{eq:inputDynamics}, has a path from an input vertex to it.

\begin{algorithm}[!htb]
    \caption{\footnotesize Distributed algorithm to verify condition $(4i)$ of Theorem~\ref{thm:struct-cont-cond}.}
    \label{alg:fdic-i}
    \footnotesize
    \begin{algorithmic}[1]
        \Procedure{Reached}{$\bA_i,\bB_i,\bE_{i,k} \neq 0,\bE_{k,i} \neq 0,r$}
            \State $\texttt{nghI}(i) \gets \{ j : \bE_{i,j} \neq 0 \}$
            \State $\texttt{nghO}(i) \gets \{ j : \bE_{j,i} \neq 0 \}$
            \State $\texttt{nghbs}(i) \gets \texttt{nghI}(i) \cup \texttt{nghO}(i)$
            \State $N_i \gets \#\{\text{SCCs of }\cD(\bA_i)\}$
            \State $\texttt{SCCs}(i) \gets \{(i,N_i)\}$
            \LComment{the subsystems interact with each other to learn how many SCCs each subsystem has, in order to find the necessary number of communication steps}
            \For{$k=1,\dots,r$}
                \label{alg:fdic-i:scc-begin}
                \ForAll{$j\in\texttt{nghbs}(i)$}
                    \State $\Call{Send}{\texttt{SCCs}(i),j}$
                    \State $\texttt{SCCs}(j) \gets \Call{Rcv}{j}$
                    \State $\texttt{SCCs}(i) \gets \texttt{SCCs}(i)\cup\texttt{SCCs}(j)$
                \EndFor
            \EndFor
            \label{alg:fdic-i:scc-end}
            \State $N \gets \sum_{j=1}^r \texttt{SCCs}(j)$
            \label{alg:fdic-i:N}
            \State $\texttt{rchd}(i) \gets \{\}$
            \Comment initialized the list of input-reached vertices
            \LComment{add the vertices with incoming edges from input vertices}
            \For{$j=1,\dots,n_i$}
                \If{$\exists k : (\bB_i)_{j,k} = 1$}
                    \State $\Call{AddTo}{x_j,\texttt{rchd}(i)}$
                \EndIf
            \EndFor
            \For{$k=1,\dots, N$}
                \label{alg:fdic-i:glob-begin}
                \wComment{transmit to outgoing neighbors which vertices that interact with them have been input reached}
                \ForAll{$j\in\texttt{nghO}(i)$}
                    \label{alg:fdic-i:com-send-begin}
                    \wComment{$M_{\bullet,l}$ is the $l$--th column of $M$}
                    \State $\Call{Send}{\{x_l: (i,x_l)\in\texttt{rchd}(i)\text{~and~}(\bE_{j,i})_{\bullet,l} \neq 0\},j}$
                \EndFor
                \label{alg:fdic-i:com-send-end}
                \LComment{~add vertices reached from the neighbors' input reached vertices}
                \ForAll{$j\in\texttt{nghI}(i)$}
                    \label{alg:fdic-i:com-rcv-begin}
                    \State $\texttt{avail}(j) \gets \Call{Rcv}{j}$
                    \State $\texttt{rchd}(i) \gets \texttt{rchd}(i) \cup \{x_t : (\bE_{i,j})_{t,l} = 1, x_l\in\texttt{avail}(j)\}$
                \EndFor
                \label{alg:fdic-i:com-rcv-end}
                \LComment{verify which vertices are reachable from the inputs by using $k$ edges between subsystems steps}
                \For{$l=1,\dots,n_i$}
                    \label{alg:fdic-i:propg-begin}
                    \State $\texttt{rchd}(i) \gets \texttt{rchd}(i) \cup \{x_t : (\bA_i)_{t,s} = 1, x_s\in\texttt{rchd}(i)\}$
                \EndFor
                \label{alg:fdic-i:propg-end}
                \label{alg:fdic-i:glob-end}
            \EndFor
            \LComment{return $\textbf{True}$ if every state vertex the $i$--th subsystem digraph is input-reached, and $\textbf{False}$ otherwise}
            \State \Return $[\#\texttt{rchd}(i) == n_i]$
        \EndProcedure
    \end{algorithmic}
\end{algorithm}

\begin{theorem}
    Algorithm~\ref{alg:fdic-i} is correct (i.e.,~it returns $\textbf{True}$ if and only if every state vertex in the $i$--th subsystem digraph is input-reached).  Further, Algorithm~\ref{alg:fdic-i} has complexity $$\cO\left(\max\left\{ r^2 , N r, N\max\limits_{i=1,\dots,r}n_i\right\}\right),$$ where $n_i$ is the dimension of the state space of the $i$--th subsystem, and $N = \sum_{i=1}^r n_i$, where $k_i$ is the number of SCCs in the $i$--th subsystem digraph.

    \vspace{-.1cm}\hfill$\diamond$
    \label{thm:fdic-a}
\end{theorem}

\begin{IEEEproof}
    Note that, since each subsystem can establish two-way communication with its neighbors, the communication graph is strongly connected, and thus the instructions in steps~\ref{alg:fdic-i:scc-begin}--\ref{alg:fdic-i:scc-end} only need to be executed (at most) $r$ times in order to receive all pairs $(\textit{id},\text{\#SCCs})$ in the system.  Subsequently the total number of SCCs, $N$, can be computed in step~\ref{alg:fdic-i:N}.

    Now, assume that each subsystem has a strongly connected state digraph $\cD(\bA_i)$.  Then, if the system has an input vertex, i.e., if $\bB_i\neq 0$, each of the state vertices of the $i$--th subsystem is added to $\texttt{rchd}(i)$ in the first iteration of the for-loop in steps~\ref{alg:fdic-i:glob-begin}--\ref{alg:fdic-i:glob-end}, namely in the for-loop~\ref{alg:fdic-i:propg-begin}--\ref{alg:fdic-i:propg-end}.  Further, note that in the case where each subsystem has a strongly connected system digraph, $N=r$, and a path from an input vertex to a state vertex contains at most $r$ edges between different subsystem digraphs.  Therefore, in this case, after $N$ iterations of steps~\ref{alg:fdic-i:glob-begin}--\ref{alg:fdic-i:glob-end} all vertices that may be reached by a path from an input vertex have been added to $\texttt{rchd}(i)$.

    Alternatively, if the $i$--th subsystem is not strongly connected, then, assume, without loss of generality, that $\bA_i$ is a block matrix, with submatrices $\bA_i^1,\dots,\bA_i^l$ along the diagonal so that $\cD(\bA_i^1),\dots,\cD(\bA_i^l)$ are strongly connected.  Further, let $\bB_i^1,\dots,\bB_i^l$ be the restriction of $\bB_i$ to the rows in used by $\bA_i^1,\dots,\bA_i^l$, respectively.  Then, consider the interconnected dynamical comprising comprising, instead of the $i$--th subsystem $(\bA_i,\bB_i)$, the subsystems $(\bA_i^1,\bB_i^1),\dots,(\bA_i^l,\bB_i^l)$ connected amongst them and to other subsystems according to $\bA_i$.  By applying this procedure to every subsystem whose state digraph is not strongly connected, we obtain an interconnected dynamical system, where each subsystem has a strongly connected digraph.  Note also, that we did not change the state digraph of the overall system, thus, a state vertex in the overall system digraph is input-reached if and only if it was input-reached in the original system digraph.  Now, since this system has $N = \#\{\text{SCCs in all subsystems of the original system digraph}\}$ subsystems, from the previous paragraph we conclude that after $N$ iterations of steps~\ref{alg:fdic-i:glob-begin}--\ref{alg:fdic-i:glob-end}, every state vertex in $\cD(\bA_i)$ that is input-reached in $\cD(\bA,\bB)$ has been added to $\texttt{rchd}(i)$.

    Thus, we have proven that for any interconnected dynamical system, every state vertex of $\cD(\bA_i)$ has a path from some input vertex in the overall system if and only if $\#\texttt{rchd}(i) = n_i$.

    Finally, we analyze the complexity of Algorithm~\ref{alg:fdic-i}: the SCCs of $\cD(\bA_i)$ can be computed in $\cO\left(n_i\right)$, and each of the steps in the for-loop~\ref{alg:fdic-i:scc-begin}--\ref{alg:fdic-i:scc-end} can be executed in constant time, which implies that the for-loop incurs in complexity $\cO\left(r\#\texttt{nghbs}(i)\right)$ which is bounded by $\cO\left(r^2\right)$.  Further, the steps~\ref{alg:fdic-i:com-send-begin}--\ref{alg:fdic-i:com-send-end} and~\ref{alg:fdic-i:com-rcv-begin}--\ref{alg:fdic-i:com-rcv-end} can be executed in constant complexity, thus these loops incur in complexity $\cO\left(\#\texttt{nghO}(i)\right)$ and $\cO\left(\#\texttt{nghI}(i)\right)$, respectively.  Finally, the for-loop in steps~\mbox{\ref{alg:fdic-i:propg-begin}--\ref{alg:fdic-i:propg-end}}, incurrs in linear complexity (on the number, $n_i$, of state variables).  So in conclusion, the complexity of Algorithm~\ref{alg:fdic-i} becomes $$\cO\left(\max\left\{ r^2 , N r, N\max\limits_{i=1,\dots,r}n_i\right\}\right).$$
\end{IEEEproof}

Next, we present a distributed algorithm to verify structural controllability when the subsystems only have access to neighboring subsystems.  Briefly, the algorithm verifies both conditions $(4i)$ and $(4ii)$ of Theorem~\ref{thm:struct-cont-cond} in a distributed manner:  on one hand, the condition $(4i)$ of Theorem~\ref{thm:struct-cont-cond} can be verified by applying Algorithm~\ref{alg:fdic-i}.  On the other hand, Theorem~\ref{thm:struct-cont-cond}--$(4ii)$ requires one to compute a maximum matching in a distributed manner.  This can be achieved by reducing the problem of finding a maximum matching to that of computing a maximum flow~\cite{ahuja_network_1993}.  However, since we only need to detect the existence of right-unmatched vertices, we only need to compute a maximum \emph{preflow}, which corresponds to a flow, where the flow on the incoming edges need not be equal to the flow on the outgoing edges of each vertex.  So,  we employ the distributed algorithm provided in~\cite{maxflow}.  In order to achieve this reduction, one takes the overall system bipartite graph and provides an orientation to each edge, from left-vertex to right-vertex. Then,  adds two extra vertices, called \emph{source} and \emph{sink}.  Finally, one adds an edge from the source to each of the left-vertices of the bipartite graph, and from each of the right-vertices to the sink and assigns to each vertex a \emph{capacity} of $1$~\cite{ahuja_network_1993}.  The computation of the maximum flow is then done distributedly, where each subsystem works to maximize the flow from the source to the sink within a \emph{region} of the graph comprising the subsystems bipartite graph, the source and the sink (note that the source and sink lie in all regions, which does not impair the distribution of the algorithm, since the systems need not keep track of the excess on the source or the sink), and any vertices in other subsystems to which the system is connected.  This is achieved through a \emph{push-relabel} algorithm, briefly described as follows: each of the vertices in a region keeps track of an \emph{excess} (which corresponds to the difference between the incoming and outgoing flow), and a \emph{label} or \emph{height}. The excess is then pushed from higher labels to lower labels increasing the flow through the edges between them until it reaches the sink, or the boundary.  Once this is achieved, the excess accumulated in the boundary is passed to the corresponding neighboring region, and the iterations begin again. However, the existence of boundary vertices limits the parallelization, as two instantiations of the algorithm can only (in general) be computed simultaneously, if the regions do not share vertices other than the source or the sink.

From this point onwards we refer to the individual instances of the \emph{parallell region discharge} algorithm presented in~\cite{maxflow}, by $\textsc{PRD}$. Further, we assume $\textsc{PRD}$ takes as parameters the digraph on which it operates, the capacity function, and the neighbors with which it shares vertices other than the source or the sink, and returns a maximal preflow on the digraph.

\begin{algorithm}[!htb]
    \caption{\footnotesize Distributed algorithm to verify condition $(4)$ of Theorem~\ref{thm:struct-cont-cond}.}
    \label{alg:fdic}
    \footnotesize
    \begin{algorithmic}[1]
        \Procedure{Controlled}{$\bA_i,\bB_i,\bE_{i,k} \neq 0,\bE_{k,i} \neq 0,r$}
            \State $\texttt{nghI}(i) \gets \{ j : \bE_{i,j} \neq 0 \}$
            \State $\texttt{nghO}(i) \gets \{ j : \bE_{j,i} \neq 0 \}$
            \State $\texttt{nghbs}(i) \gets \texttt{nghI}(i) \cup \texttt{nghO}(i)$
            \LComment{verify if every state vertex is input-reached by deploying Algorithm~\ref{alg:fdic-i}}
            \State $\texttt{rchd}(i) \gets \Call{Reached}{\bA_i,\bB_i,\bE_{i,k} \neq 0,\bE_{k,i} \neq 0,r}$
            \label{alg:fdic:reached}
            \LComment{we set up the graph for applying the Parallel region discharge, where $s$ and $t$ correspond to the source and sink, respectively, the $x$ and $u$ vertices correspond to state and input vertices.  The upper index $i$ is the index of the subsystem they belong to, and the upper index $R$ and $L$ indicate if they are right or left vertices}
            \State $\cV_i \gets \{ s, t\} \cup \{ x^{i,L}_k , x^{i,R}_k\}_{k=1}^{n_i} \cup \{ u^{i}_k\}_{k=1}^{m_i}$
            \label{alg:fdic:graphsetup-begin}
            \State $\cE_{i,i} \gets \{ (x^{i,L}_j,x^{i,R}_{j'}) : (\bA_i)_{j',j} = 1 \} \cup \{ (u^{i}_j,x^{i,R}_{j'}) : (\bB_i)_{j',j} = 1 \}$
            \ForAll{$ j \in \texttt{nghbs}(i)$}
                \State $\cV_{i,j} \gets \{  x^{j,R}_k : (\bE_{j,i})_{k,\bullet} \neq 0\}$
                \State $\cE_{i,j} \gets \{ (x^{i,L}_l,x^{j,R}_k) : (\bE_{j,i})_{k,l} = 1 \}$
                \State $\cV_{j,i} \gets \{  x^{j,L}_k : (\bE_{j,i})_{\bullet,k} \neq 0\}$
                \State $\cE_{j,i} \gets \{ (x^{j,L}_k,x^{i,R}_l) : (\bE_{i,j})_{l,k} = 1 \}$
            \EndFor
            \State $\cE_{s} \gets \{s\} \times \{ x^{i,L}_k \}_{k=1}^{n_i}$
            \State $\cE_{t} \gets \{ x^{i,R}_k\}_{k=1}^{n_i} \times \{t\}$
            \State $\cE \gets \cE_s\cup\cE_t\cup\cE_{i,i}\cup\bigcup\limits_{j\in\texttt{nghbs}(i)}(\cE_{j,i}\cup\cE_{j,i})$
            \State $\cV \gets \cV_i \cup \bigcup\limits_{j\in\texttt{nghbs}(i)}(\cV_{j,i}\cup\cV_{i,j})$
            \State $\cD \gets (\cV,\cE)$
            \Function{$c$}{$e \in \cE_i$}
                \State \Return $1$
                \Comment{all edges have unitary capacity}
            \EndFunction
            \label{alg:fdic:graphsetup-end}
            \LComment{Deploy a Parallel Region Discharge algorithm to obtain a preflow $f$ on $\cD$, with capacity function $c$}
            \State $f \gets \textsc{PRD}(\cD,c,\texttt{nghbs(i)})$
            \label{alg:fdic:prd}
            \State $\texttt{mchd}(i) \gets [\sum\limits_{e\in\cE_t} f(e) ==  n_i]$
            \label{alg:fdic:right-unmatch}
            \LComment{check if the whole system can be structurally controllable according to whether the conditions are satisfied in the current system or not}
            \State $\texttt{ctld}(i) \gets \texttt{rchd}(i)\wedge\texttt{mchd}(i)$
            \label{alg:fdic:init-guess}
            \For{$k=1,\dots,r$}
            \label{alg:fdic:un-begin}
                \LComment{reconsider the controllability of the overall system, in light of the data received from the neighbors}
                \ForAll{$j\in\texttt{nghbs}(i)$}
                    \State $\Call{Send}{\texttt{ctld}(i),j}$
                    \State $\texttt{ctld}(j) \gets \Call{Rcv}{j}$
                    \State $\texttt{ctld}(i) \gets \texttt{ctld}(i)\wedge\texttt{ctld}(j)$
                \EndFor
            \EndFor
            \label{alg:fdic:un-end}
            \State \Return $\texttt{ctld}(i)$
            \LComment{return \textbf{True} if the system is structurally controllable and \textbf{False} otherwise}
        \EndProcedure
    \end{algorithmic}
\end{algorithm}

\begin{theorem}
    Algorithm~\ref{alg:fdic} is correct, i.e.,~it verifies  Theorem~\ref{thm:struct-cont-cond}--$(4)$ in a distributed fashion.  Further, it has a computational complexity of
    $$\cO\left(\max \{ r^2 , Nr , N \max\limits_{i=1,\dots,r} n_i, r\beta^2 \max\limits_{i = 1,\dots,r} n_i^3\}\right),$$
    where $\beta$ is the number of boundary vertices, and the remaining variables are the same as described in Theorem~\ref{thm:fdic-a}.
    \hfill$\diamond$
    \label{thm:algc-fdic}
\end{theorem}

\begin{IEEEproof}
    In order to verify the correctness of Algorithm~\ref{alg:fdic}, we have to check if both conditions $(4i)$ and $(4ii)$ of Theorem~\ref{thm:struct-cont-cond} are verified. Furthermore, in order to perform this verification in a distributed manner, each subsystem must verify that all vertices in its digraph are input-reached in $\cD(\bA,\bB)$, which is done by employing Algorithm~\ref{alg:fdic-i} in step~\ref{alg:fdic:reached}; and that none of its state vertices are right-unmatched w.r.t.\ some maximum matching of $\cB(\bA,\bB)$. In addition, recall that it was already argued in the proof of Theorem~\ref{thm:sufP1alg} that the for-loop in steps~\ref{alg:fdic:un-begin}--\ref{alg:fdic:un-end} determines if these conditions are violated in any of the subsystems.

    Now, to verify that Algorithm~\ref{alg:fdic} determines if there are right-unmatched vertices in the $i$--th subsystem, we note that in steps~\ref{alg:fdic:graphsetup-begin}--\ref{alg:fdic:graphsetup-end} we generate the digraph $\cD$ comprising the right- and left-vertices of the $i$--th subsystem bipartite graph, and the boundary vertices of the $i$--th region according to the precepts of~\cite{maxflow}.  Once the digraph $\cD$ is computed, we apply $\textsc{PRD}$ to it in step~\ref{alg:fdic:prd}, thus obtaining a preflow from source to sink on $\cD$ which is maximum amongst preflows on the whole graph.  By the guarantees provided in~\cite{maxflow}, together with the equivalence between the maximum matching and maximum flow problems, presented in~\cite{ahuja_network_1993} we guarantee that $\sum\limits_{e\in\cE_t}f(e)$ is equal to the number of right-matched vertices in a maximum matching of the system bipartite graph, that are state vertices of the $i$--th subsystem.  So, by comparing $\sum\limits_{e\in\cE_t}f(e)$ with $n_i$ in step~\ref{alg:fdic:right-unmatch}, we are able to infer if there are right-unmatched vertices in the $i$--th subsystem w.r.t.~some maximum matching $\cB(\bA,\bB)$.  Thus the algorithm returns $\textbf{True}$ if and only if every state vertex of the system digraph is input-reached, and there are no right-unmatched vertices in the system bipartite graph.

    Now, since all of the steps of the algorithm have linear complexity except for step~\ref{alg:fdic:reached} and step~\ref{alg:fdic:prd}, the complexity of Algorithm~\ref{alg:fdic} is given by the maximum of these.
    Knowing that step~\ref{alg:fdic:reached} has a complexity of $\cO\left(\max\left\{ r^2 , N r, N\max\limits_{i=1,\dots,r}n_i\right\}\right)$ (see Theorem~\ref{thm:algc-fdic}, where $N$ is the number of SCCs on each of the subsystems), all that remains to infer is the complexity of step~\ref{alg:fdic:prd}.  This algorithm, as described in~\cite{maxflow},   has complexity $\cO\left(n_i^3\right)$~\cite{ahuja_computational_1997,goldberg_partial_2008}, and the necessary iterations of the region discharge that each subsystem has to complete, can be bounded by $\beta^2$, where $\beta$ is the number of boundary vertices in the whole system bipartite graph (that is, the number of vertices in the bipartite graph that have to be shared by several subsystems).  Also, in the worst-case scenario, where each of the subsystems is connected to every other subsystem, the region discharge steps have to be executed sequentially. So, the complexity of step~\ref{alg:fdic:prd} is given by $\cO\left(r\beta^2 \max\limits_{i = 1,\dots,r} n_i^3\right)$, resulting in an overall complexity of $\cO\left(\max\left\{ r^2 , N r, N\max\limits_{i=1,\dots,r}\left\{n_i , r\beta^2 n_i^3 \right\}\right\}\right)$.
\end{IEEEproof}

\section{Illustrative Examples}\label{illustrativeexample}

\begin{figure}[htb]   
    \centering
\includegraphics[scale=0.55]{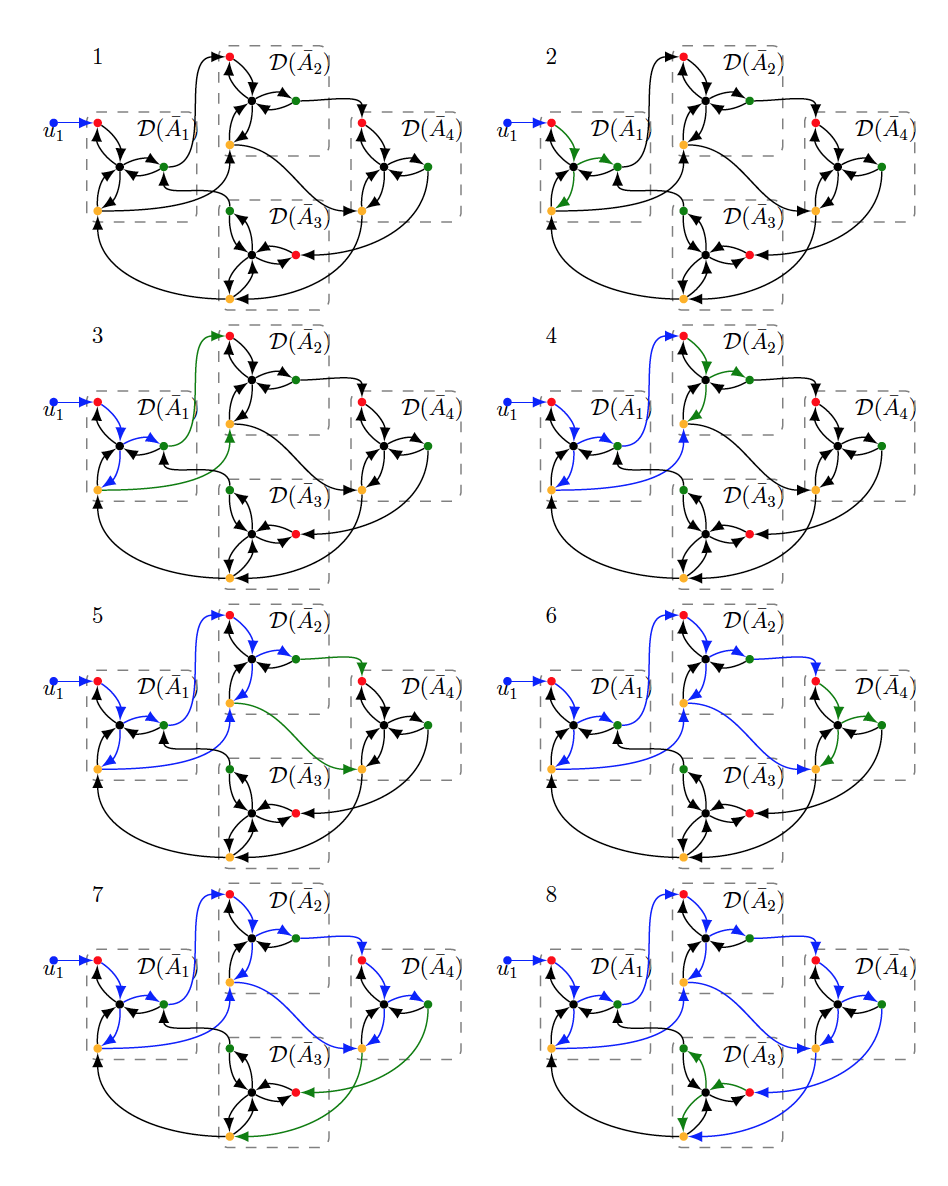}
    \caption{\footnotesize Example of the procedure of Algorithm~\ref{alg:fdic-i} applied to the system digraph presented in Figure~\ref{fig:reachability}--1, comprising $4$ different subsystems (depicted inside of the dashed boxes), only one of which has an input edge (labeled $u_1$).  In each subfigure, the blue edges represent those that comprised a path from an input vertex, and the green edges denote those that were added in this iteration or communication step of the algorithm.  Finally, the even-labeled subfigures correspond to an iteration of Algorithm~\ref{alg:fdic-i}, and the odd-labeled ones correspond to a communication step between subsystems.}
    \label{fig:reachability}
\end{figure}

Let us consider the control system, whose system digraph is depicted in Figure~\ref{fig:reachability}--$1$. This consists of four interconnected
subsystems, whose state digraphs are enclosed by grey dashed boxes. We would like to assess if this system is structurally controllable from the single control input $u_1$, considering only the locally available information, i.e., in a distributed fashion. To solve this problem, we apply  Algorithm~\ref{alg:fdic} (and Algorithm~\ref{alg:fdic-i}, which is required as a subroutine) that is executed on subsystem level, hence, emphasizing its distributed nature.

In Figure~\ref{fig:reachability}--1, we present the digraph associated to an interconnected dynamical system comprising four subsystems.  Since only subsystem $(\bA_1,\bB_1)$ has an input vertex, it readily follows from Theorem~\ref{thm:struct-cont-cond} that none of the other subsystems can be structurally controllable.  Now, we deploy Algorithm~\ref{alg:fdic} to verify the structural controllability of the interconnected dynamical system.  After the initialization steps are completed, we deploy Algorithm~\ref{alg:fdic-i}, the iterations of which can be seen in Figure~\ref{fig:reachability}--2 to Figure~\ref{fig:reachability}--8: in each iteration (even-labeled subfigures) new vertices are seen to be input-reached (the targets of the green edges), and in each communication step (odd-labeled subfigures) the subsystems interact to its outgoing neighbors which of their vertices are reached after the iteration has been completed.  As can be seen in Figure~\ref{fig:reachability}--8, all vertices have been reached after four iterations, which in this case, corresponds to the number of SCCs in all subsystems.

\begin{figure}[htb]    \centering
    \centering
\includegraphics[scale=0.45]{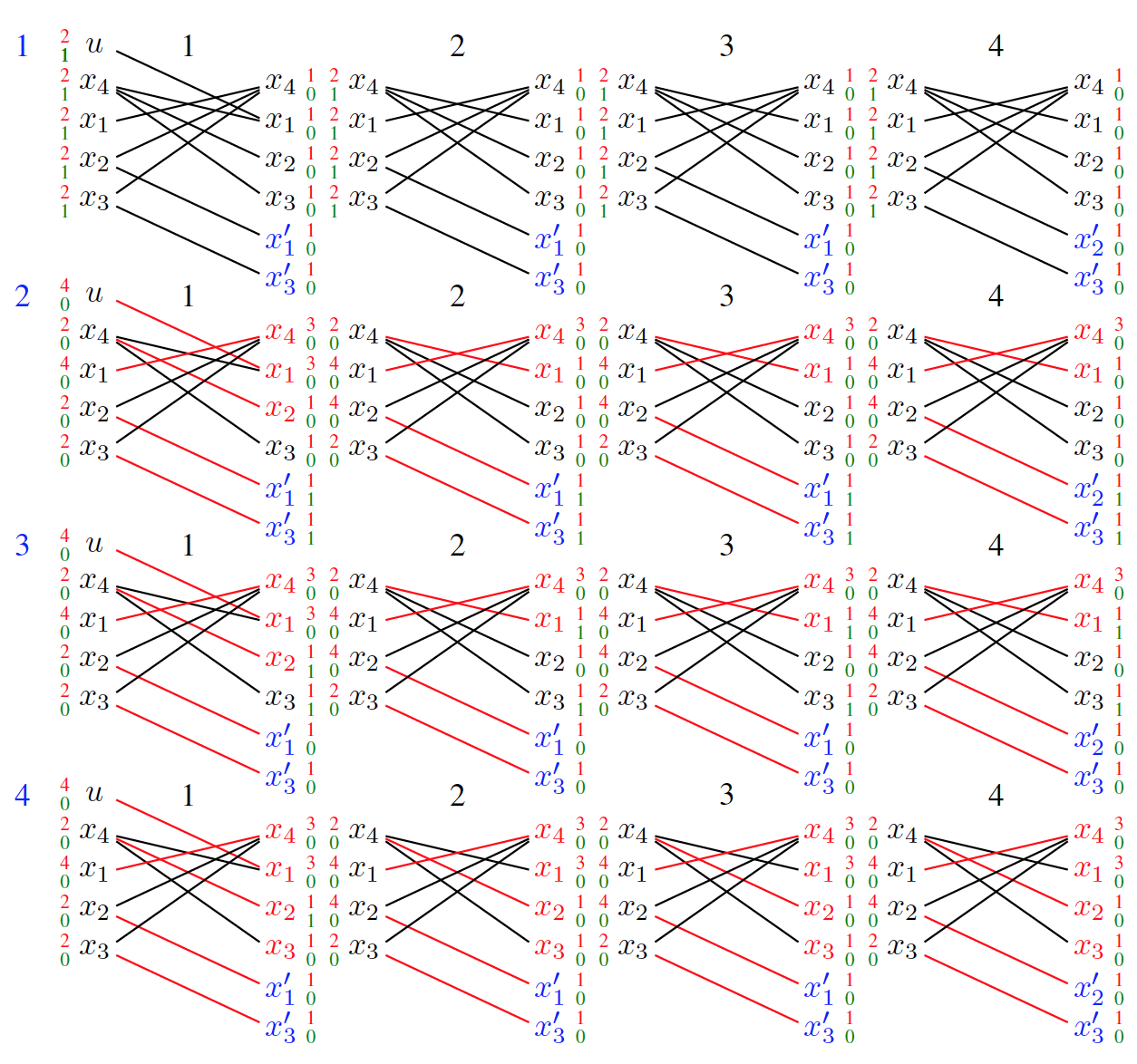}
    \caption[]{\footnotesize Example of the $\textsc{PRD}$ algorithm applied to the verification of structural controllability of the system presented in Figure~\ref{fig:reachability}--1.  For convenience of referencing, the vertices are given labels rather than colors ($x_4$ being the black vertex in each subsystem and the others can be easily inferred).  At the left of each left-vertex and at the right of each right-vertex we insert two numbers, the one in green represents the excess of the corresponding vertex, whereas the red number represents its label.  Edges in red represent those where the capacity has been saturated; and right-vertices in red represent the ones for which the edge to the sink has been saturated.  Finally, the vertices in blue represent vertices that belong to other regions, i.e., boundary vertices.  In order to simplify, in this, we do not include boundary vertices from incoming subsystems. Note also, that in this instance, we can run the algorithm in all regions simultaneously, since by not considering the incoming edges from other regions in the region graph, we do not allow for flow to be sent back through these edges.}
    \label{fig:reg-disch}
\end{figure}

In Figure~\ref{fig:reg-disch}, we consider an example of a run of the region discharge algorithm running on the bipartite graph associated to the digraph in Figure~\ref{fig:reachability}--1.  In this example, we begin applying $\textsc{PRD}$ in step $1$ by initializing the labels at $2$ for each left-vertex, and at $1$ for each right-vertex; we also saturate all edges from the source, which makes it so that all of the left-vertices start with an excess of $1$.  By successive pushing and relabeling, we reach the configuration in $2$ where all the excess has either been pushed to the sink (and thus the corresponding right-vertex is presented in red) or to the boundary of the region.  In step $3$, we discharge the excess from the boundary into the adjacent regions so that, for example, the right-vertex $x_3$ in each of the regions has now an excess of $1$.  Finally, by applying push-relabel again in each of the regions, we reach step $4$ where all the edges from right-vertices to the sink have been saturated (and are thus displayed in red) showing that there is a maximum preflow saturating all edges to the sink, and equivalently that there is a maximum matching with no right-unmatched vertices.  So, in combination with the analysis of Figure~\ref{fig:reachability} we conclude that the interconnected dynamical associated to the digraph system in Figure~\ref{fig:reachability}--$1$ is structurally controllable.

\section{Conclusions and Further Research}
\label{conclusions}

In this paper, we have provided several necessary and/or sufficient conditions to verify structural controllability for interconnected linear time-invariant dynamical systems based on the local information accessible to each subsystem.  Subsequently, we have provided distributed and efficient (i.e., polynomial in the dimension of the state and input) algorithms to verify a necessary and sufficient condition for structural controllability.  The results presented readily extend to discrete time-invariant interconnected dynamical systems, since the controllability criterion stays the same.  Further, by duality between controllability and observability the results also apply to structural observability verification of discrete/continuous linear time-invariant interconnected dynamical systems.  Whereas the results presented pertain to verifying structural conditions, it would be of interest to address design problems; for instance, which state variables need to be actuated, or which inputs should be used, to ensure a given structural property.  On the other hand, it would be of interest to understand if the conditions provided could be adapted to the case where some of the entries in the structure of the subsystems and their interconnections are known exactly (which corresponds to the case where only some of the components of the overall system are assumed to be reliable).  Ultimately, such an extension would shed light on the relationship between structural and non-structural system-theoretic properties; hence, leading to a better understanding of the resilience and performance of interconnected dynamical systems.

\section{Acknowledgements}

The first author would like to thank ACCESS Linnaeus Center, for their hospitality, as much of the writing of this paper and the research that led to it were carried out while visiting KTH in the Spring of 2014.

\bibliographystyle{IEEEtran}
\bibliography{IEEEabrv,bibliography}

\clearpage

\end{document}